\documentclass[11pt]{article}

\usepackage{amsmath, amssymb}
\usepackage{graphicx}

\newcommand{\cqfd}{{\nobreak\hfil\penalty50\hskip2em\hbox{}\nobreak\hfil
$\square$\qquad\parfillskip=0pt\finalhyphendemerits=0\par\medskip}}

\newcommand{\R}{\mathbb{R}}
\newcommand{\Z}{\mathbb{Z}}

\newcommand{\dt}{\partial_t}

\newcommand{\dx}{\partial_x}

\newcommand{\ds}{\displaystyle}

\newcommand{\eps}{ \varepsilon}

\newcommand{\pr}{{\bf \textit{Proof: }}}

\newtheorem{theorem}{Theorem}[section]

\newtheorem{proposition}{Proposition}[section]
 \newtheorem{lemma}{Lemma}[section]
\newtheorem{definition}{Definition}[section]
\newtheorem{remark}{Remark}[section]

\title{ \bf 
 Strong Stability with respect to weak limit 
  for  a Hyperbolic System  
 arising from   Gas Chromatography}

\author {
C. Bourdarias
\thanks{Universit\'{e} de Savoie, LAMA, UMR CNRS 5127,
 73376 Le Bourget-du-Lac.
bourdarias@univ-savoie.fr},
M. Gisclon
\thanks{Universit\'{e} de Savoie, LAMA, UMR CNRS 5127,
73376 Le Bourget-du-Lac.
 gisclon@univ-savoie.fr }
and S. Junca
\thanks{IUFM et Universit\'{e} de Nice,
 Labo. JAD, UMR CNRS 6621, Parc Valrose, 06108, Nice.
 junca@unice.fr}
}

\date{}

\begin{document}

\bibliographystyle{plain}

\maketitle

\abstract{ 
We investigate a system related to a particular isothermal gas-solid
chromatography process, called ``Pressure Swing Adsorption'', with two species
and instantaneous exchange kinetics. This system presents the particularity to
have a  linearly degenerate eigenvalue: this allows the velocity of the
gaseous mixture to propagate high frequency waves. In case of smooth
concentrations with a general isotherm, we prove $L^1$ stability for
concentrations with respect to weak limits of the inlet boundary velocity.
Using the Front Tracking Algorithm (FTA), we prove a similar
result for
concentrations with bounded variation
(BV) under some convex assumptions on the isotherms.
In both cases  we show that  high frequency oscillations with large amplitude of
the inlet velocity can propagate without affecting the concentrations.
}




\section{Introduction}

``Pressure Swing Adsorption (PSA) is a technology  used to separate some species from a gas under
pressure according to the  molecular characteristics and affinity of the species for an adsorbent
material. Special adsorptive materials (e.g. zeolites) are used as a molecular sieve, preferentially
adsorbing the undesired gases at high pressure. The process then swings to low pressure to desorb
the adsorbent material'' (source: Wikipedia).

A typical PSA system involves a cyclic process where a number of connected vessels containing
adsorbent material undergo successive pressurization and depressurization steps in order to produce
a continuous stream of purified product gas. We focus here on a step of the cyclic process,
restricted to isothermal behavior.\\ 
As in general fixed bed chromatography, each of the $d$ species ($d\geq 2$) simultaneously exists
under two phases, a gaseous and movable one with velocity $u(t,x)$ and concentration $c_i (t,x)$ 
or a solid (adsorbed) other with concentration $q_i (t,x)$, $1\leq i\leq d$. We assume that mass 
exchanges between the mobile and the stationary phases are infinitely fast, thus the two phases are
constantly at  composition e\-qui\-li\-brium: the concentrations in the solid  phase are given by
some relations  $q_i=q_i^*(c_1,...,c_d)$ where the functions $q_i^* $ are the so-called
e\-qui\-li\-brium isotherms. A theoretical study  of a model with finite
exchange kinetics was presented in \cite{B92} and a numerical approach was developed in \cite{B98}. 

\newpage

In gas chromatography, velocity variations accompany changes in gas composition, especially in the
case of high concentration solute: it is known as the sorption effect. In the present model, the
sorption effect is taken into account through a constraint on the pressure (or on the density in
this i\-so\-ther\-mal case). See \cite {RAA70} and \cite{Ru84} for a precise description of the
process and \cite{BGJ08} for a survey on various related models. 

The system for two species ($d=2$) with three unknowns $(u,c_1,c_2)$ is:
\begin{eqnarray}
\dt (c_1+q^*_1(c_1,c_2))+\dx(u\,c_1)&=&0, \label{un}\\
\dt (c_2+q^*_2(c_1,c_2))+\dx(u\,c_2)&=&0, \label{deux}\\
c_1+c_2&=& \rho(t), \label{trois}
\end{eqnarray}
with suitable initial and boundary data.
The function $\rho$ represents the \textit{given} total density of the mixture. The experimental
device is realized so that it is a given function depending only upon time and in the sequel we
assume that $\rho\equiv 1$ (which is not really restrictive from a theoretical point of view).
First existence results of large solutions satisfying some entropy criterium in the case of two
chemical species were obtained in \cite{BGJ06,BGJ07}.
In the previous system, it appears that we can expect strong singularities  with respect to time for
the velocity $u$. For instance,  let $c_1(t,x)\equiv \underline{c}_1$ be a  constant,
$c_2(t,x)\equiv 1-\underline{c}_1$, $u(t,x)\equiv u_b(t)$ where $u_b$ is any $L^\infty$
function, then $(c_1,c_2,u)$ is a weak solution of (\ref{un}),(\ref{deux}),(\ref{trois}).
So we can build solutions with a strong oscillating velocity for this system. Furthermore
high oscillations of the  incoming velocity $u_b$ slightly perturb the concentration as we will
see. Notice that we seek  positive solutions $(c_1,c_2)$, thus,  in view of (\ref{trois}) with
$\rho\,\equiv 1$, $c_1$, $c_2$ must satisfy $$0\leq c_1,\,c_2\leq 1.$$ We use the following
notations, introduced in \cite{BGJ07}: we set $c=c_1 \in [0,1]$  and
\begin{eqnarray*}
q_i(c)&=&q_i^ *(c,1-c),\quad i=1,2,\\
h(c)&=&q_1(c)+q_2(c),\\
I(c)&=&c+q_1(c).
\end{eqnarray*}
Adding (\ref{un}) and (\ref{deux}) we get, thanks to (\ref{trois}):
$$\dt (q_1(c)+q_2(c))+\dx u=0,$$
 thus our purpose is to study the following system:
\begin{equation}\label{sysad}
\left\{   \begin{array}{ccl}  \vspace{2mm}
\dt I(c)+  \dx  (u\, c) & = & 0,\\
\dt h(c) + \dx u & = & 0, 
\end{array}\right.
\end{equation}
supplemented by initial and boundary data:
\begin{equation} \label{sysad0}
\left\{   \begin{array}{ccl}
 \vspace{2mm}c(0,x)&=&c_0(x) \in [0,1],  \quad x > 0,\\
 \vspace{2mm}c(t,0) &=&c_b(t)   \in [0,1],\quad t>0,\\
 u(t,0)&=&u_b(t)>0,\quad t>0.
\end{array}\right.
\end{equation}
Notice that we assume in (\ref{sysad0}) an incoming flux at the boundary, i.e. $\forall t>0,\
u_b(t)>0$. In the case where the first species is inert, that is $q_1=0$, the $I$ function reduces
to identity.

System (\ref{sysad}) has  a null eigenvalue as the system exposed in \cite{BJ97}, but,
instead \cite{BJ97}, we cannot reduce this system to  a single equation for general solutions with
shocks. In \cite{BR03} is studied another interesting 2x2 system with a linearly degenerate
eigenvalue which modelises some traffic flow. As in \cite{CG99,CG01,CGM04,CGM03,M04}, the zero
eigenvalue makes possible the existence of stratified solutions or the propagation of
large-amplitude high frequency waves. Usually, for genuinely nonlinear conservation laws, only high
oscillating solutions with small amplitude can propagate: see for instance \cite{DM85,CJR06}.

In this paper we prove,  for large data, that the velocity is a stratified solution in the following
sense: $u(t,x)=u_b(t)\, v(t,x)$ where $v$ is as regular as the concentration $c$ and more than  the
boundary data $u_b$. This decomposition for the velocity allows high oscillations with large
amplitude for velocity to propagate, without affecting the concentration. For this quasilinear
system we have propagation of high oscillations with large amplitude for velocity as in a semilinear
system, see for instance \cite{JMR93,J98-1,J08}, and we have strong profile for $u$ with double
scale as in
\cite{J98-2}.

This also permits to pass to the weak limit for $u$ at the boundary and to the strong limit in the
interior for the concentration. For the smooth case, we have no restriction on the isotherms, but
for the realistic case with shock-waves, we restrict ourselves to the classical treatment   of
hyperbolic systems: eigenvalues are linearly degenerate or genuinely nonlinear. Furthermore we
obtain better interaction estimates when the shock and rarefaction  curves are monotonous. It is
the case for instance for an inert gas and an active gas with the Langmuir isotherm. We conjecture
that our result is still valid for general isotherms with piecewise genuinely nonlinear eigenvalue.
  
The paper is organized as follows. In Section \ref{sh} we recall some basics results from
\cite{BGJ07} concerning hyperbolicity, entropies, weak entropy solutions of
System (\ref{sysad}).

In Section \ref{skf}, we  study the case where concentration is smooth and the velocity is only
$L^\infty$.

In the remainder of the paper we study the case with only $BV $ concentrations. In short section
\ref{sFTA} we briefly expose the Front Tracking Algorithm (FTA) for  System
(\ref{sysad}).

Section \ref{sRP} is devoted to the study of both shock and rarefaction curves. We state 
the assumptions that we need to perform estimates with the Front Tracking Algorithm.  These 
assumptions restrict us to  convex (or concave) isotherms and we give some examples from chemistry.
We obtain the fundamental interaction estimates in Section \ref{sIE} and $BV$ estimates for $v$ in
Section \ref{suBV}. Finally, we obtain strong stability for concentration with respect to weak limit
on the boundary velocity in Section \ref{sks}.


\section{Hyperbolicity  and entropies} \label{sh}


In order the paper to be self contain, we recall without any proof some results exposed in
\cite{BGJ07}.

It is well known that it  is possible to analyze the system of Chromatography, and thus System
(\ref{sysad}), in terms of hyperbolic system of P.D.E. provided we exchange the time and space
variables and $u >0$: see \cite{RAA86} and also \cite{RSVG88} for instance. In this framework the
vector state will be 
$U=\left(
\begin{array}{l}
u \\
m
\end{array}\right)$ 
where $m=u\,c $ is the flow rate of the first species. In this vector state, $u$ must be understood
as $u\,\rho$, that is the total flow rate.

In the sequel, we will make use of the function $f=q_1\,c_2-q_2\,c_1$ introduced by Douglas and {\em
al.} in \cite{DCRBT88}, written here under the form 
\begin{eqnarray}  \label{eqf}
    f(c)& = & q_1\,c_2-q_2\,c_1= q_1(c)-c\,h(c).
\end{eqnarray} 
Any equilibrium isotherm related to a given species is always increasing with respect to the
corresponding concentration (see \cite{DCRBT88}) i.e. $\ds\frac{\partial q^*_i}{\partial c_i}\geq
0$.
Since $c = c_1$ and $c_2 = 1 - c$, it follows: 
\begin{equation}\label{qprime}
q'_1\geq   0  \geq q'_2.
\end{equation}
Let us define the function $H$ by
\begin{eqnarray}  \label{H}
H(c) & =1+(1-c)\,q'_1(c)-c\,q'_2(c)& =1+q_1'(c)-ch'(c).
\end{eqnarray}
From (\ref{qprime}), $H$ satisfies $H\geq 1$ and we have the following relation 
 between $f$, $H$ and $h$:
\begin{eqnarray*}\label{fsec}
f''(c)&=&H'(c)-h'(c).
\end{eqnarray*}

\newpage
\subsection{Hyperbolicity}


Concerning hyperbolicity, we  refer to \cite{D00,S96,Sm94}. System (\ref{sysad}) takes the form

\begin{equation}\label{sysadum}
\partial_x U +\partial_t \Phi(U)=0\hbox{  with }
U=\left(
\begin{array}{l}
u \\
m
\end{array}\right)
\hbox{  and }
\Phi(U)=\left(
\begin{array}{l}
h(m/u) \\\\
I(m/u)
\end{array}\right).
\end{equation}
The eigenvalues are: 
 \begin{eqnarray*}  
 \label{eigenvalues} 
  0  &\mbox{  and } &
\lambda=\ds\frac{H(c)}{u},
\end{eqnarray*}
 thus in view of (\ref{H}) the system is strictly hyperbolic. The zero eigenvalue is of course
linearly degenerate, moreover the right eigenvector 
$r=\left(\begin{array}{c} h'(c) \\
1+q'_1(c)
\end{array} \right )$ 
associated to $\lambda$ satisfies $\ds d\lambda \cdot
r=\ds\frac{H(c)}{u^2}\,f''(c)$, thus $\lambda$
is genuinely nonlinear in each domain where $f''\neq 0$.


\begin{proposition}[\cite{BGJ07} Riemann invariants]\label{PropRi}~\\
System (\ref{sysad}) admits the two Riemann invariants:
\begin{eqnarray*} 
c  \quad \hbox{ and } \quad w= \ln(u)+g(c)=L+g(c),
\qquad \mbox{ where }  g'(c)=\ds\frac{-h'(c)}{H(c)} \quad \mbox{and }  
L=\ln(u).
\end{eqnarray*}
Furthermore this system can be rewritten for smooth solutions as:
\begin{eqnarray} 
\label{smooth}
\partial_x c + \frac{H(c)}{u}\, \partial_t c = 0,\quad
\quad \partial_x( \ln(u)+g(c) )= \partial_x w= 0.
\end{eqnarray}
\end{proposition}


\subsection{Entropies}\label{subsecEntropies}


Dealing with entropies,  it is more convenient, as  shown in \cite{BGJ07}, to
work with the functions
\begin{eqnarray*} 
 G(c)=\exp(g(c)),    &  W=\exp(w)=u\,G(c).
\end{eqnarray*}
Notice that $G$ is a positive solution of 
 $ H G'+h'G=0$.
\\
Denote $E(c,u)$ any smooth entropy and $Q=Q(c,u)$ any associated entropy flux.
Then, for smooth solutions,
$\partial_x E +\partial_t Q=0$. Moreover:

\begin{proposition}[\cite{BGJ07} Representation of all smooth
entropies]\label{smoothentr}~\\
The smooth entropy functions for  System (\ref{sysad}) are given by
\begin{eqnarray*} 
E(c,u) & =&\phi(w)+u \,\psi(c)
\end{eqnarray*}
where $\phi$ and $\psi$ are any smooth real functions.
The corresponding entropy fluxes satisfy 
\begin{eqnarray*}
Q'(c)  & = & h'(c)\, \psi(c) +  H(c)\, \psi'(c) .
\end{eqnarray*}
\end{proposition}

Moreover, in \cite{BGJ3}, the authors looked for convex entropies for  System
(\ref{sysadum}) (i.e. System (\ref{sysad}) written in the $(u,m)$ variables) in
order to  get a kinetic formulation.
 The next proposition gives us a family of degenerate convex 
entropies independently of convexity of $f$ or of the isotherms.
\begin{proposition}[\cite{BGJ07} Existence of degenerate convex entropies]~\\
If $\psi$ is convex  or degenerate convex,  i.e. $\psi'' \geq 0$,
then $E=u\,\psi(c)$ is a degenerate convex entropy.
\end{proposition}
There are some few cases (water vapor or ammonia for instance) where the
isotherm is convex. There is also the important case with an inert carrier gas
and an active gas with a concave or convex isotherm (see
\cite{BGJ06,BGJ07,BGJ08}). In these cases, the next proposition 
ensures the  existence of  $\lambda$-Riemann invariants which  are also strictly
convex entropies. In such cases, $w$ is monotonous with respect to $x$ for any
entropy solution.
\begin{proposition}[\cite{BGJ07} When $\lambda$-Riemann invariant is a convex
entropy\label{Prop2RE}]~\\
There are strictly convex entropy of the form $E=\phi(w)$
if and only if  $G''$ does not  vanish.
\\
More precisely,  for $\alpha >0$,  $E_\alpha(c,u)= u^\alpha\, G^\alpha(c)  $
is an increasing entropy with respect to the Riemann invariant $W$. It is strictly convex for
$\alpha > 1$ if
$G'' >0$ and for $\alpha < 1$  if $G'' <0$.
\end{proposition}
Unfortunately, when  $G$ has an inflexion point such system does not admit any
strictly convex entropy. 
When one gas is inert, it is always the case   if the sign of the second
derivative of the isotherm changes.
See for instance \cite{BGJ07} for the BET isotherm.
\begin{remark}
In general, System (\ref{sysad}) is not in the Temple class. It is the case if and only if $f''$
does not vanish and $\partial_xW=0$ for  all entropy solution (\cite{BGJ5}). For
instance, System (\ref{sysad}) with two linear isotherms is in the Temple class.
\end{remark}

\begin{proposition}[\cite{BGJ07} Non Existence of  strictly convex entropy]
\label{Propconvex}~\\
If sign of $G''$ changes
then System (\ref{sysad}) does not admit  strictly  convex  smooth entropy.
\end{proposition}


\subsection{Definition of weak entropy solution}


  We have seen that there are two families of entropies:
  $u\,\psi(c)$ and $\phi(u\,G(c))$. 
\\
 The first family is degenerate convex (in variables  $(u,uc)$) 
 provided $\psi''\geq 0$. 
 So, we seek after weak entropy solutions which satisfy 
 $ \partial_x \left(u\,\psi(c) \right) + \partial_t Q(c) \leq 0$ in 
 the distribution sense.
\\
 The second family is not always convex. There are only two interesting cases,
 namely $\pm G''(c) > 0$ for all $c \in [0,1]$.
 When $G''>0$ and $\alpha > 1$, we expect to have 
 $\partial_x ( u\,G(c))^\alpha \leq 0 $ from Proposition 
\ref{Prop2RE}. 
But, the mapping $ W \mapsto W^\alpha$ is increasing on $\R^+$.
So, the last inequality reduces to 
  $\partial_x ( u\,G(c)) \leq 0 $.
\\
In the same way, if $ G''<0$, we get $ \partial_x ( u\,G(c)) \geq 0 $.
\\
Now, we can state a mathematical definition of weak entropy solutions.
\begin{definition}\label{defwes}
 Let be $T>0$, $X>0$, $u \in L^\infty((0,T)\times(0,X), \R^+)$, $ 0 \leq c(t,x) \leq \rho \equiv 1$
for almost all $(t,x) \in(0,T)\times(0,X)$. Then $(c,u)$ is a {\bf  weak entropy solution} 
 of System (\ref{sysad})-(\ref{sysad0}) with respect to the family of entropies
$u\,\psi(c)$ if, for all  convex
(or degenerate
convex) $ \psi $:
 \begin{eqnarray} \label{ineqentropies}
 \frac{\partial}{\partial x}\left(u\,\psi(c) \right) 
 +  \frac{\partial}{\partial t} Q(c) & \leq & 0,
\end{eqnarray}
in $\mathcal{D}'([0,T[\times[0,X[)$, where $Q'=H\psi'+h'\psi$, that is, for all
$\phi\in\mathcal{D}([0,T[\times[0,X[)$:
$$\int_0^X \int_0^T
\left( u\,\psi(c)\,\partial_x\phi + Q(c)\,\partial_t\phi\right) \,dt\,dx + \int_0^T
u_b(t)\,\psi(c_b(t))\,\phi(t,0)\,dt + \int_0^X Q(c_0(x))\,\phi(0,x)\,dx\geq 0.$$
\end{definition}
\begin{remark}
If $\pm G'' \geq 0$ then $u\,\psi = \pm u\,G(c)$ is a degenerate convex entropy,
with entropy flux $ Q\equiv 0$, contained in the family of entropies
$u\,\psi(c)$. So, if $G''$ keeps a constant sign on $[0,1]$, $(c,u)$ has to 
satisfy: 
  \begin{equation} \label{IRdecay}
 \pm \frac{\partial}{\partial x} \left(u\,G(c) \right)   \leq  0,
      \quad \mbox{ if } \pm G'' \geq 0 \mbox{ on $[0,1]$.}
\end{equation} 
Notice that the entropies $u\,\psi(c)$ and the entropy $u\,G(c)$ are linear with
respect to the velocity $u$.
\end{remark}


\subsection{About the Riemann Problem}


The implementation of the Front Tracking Algorithm used extensively from Section
\ref{sFTA} requires some results about the solvability of the following  Riemann problem:
\begin{eqnarray}
\left\{   \begin{array}{ccc}  \vspace{2mm}
\dx u +\dt h(c)&=&0,\\
\dx  (u c)+\dt I(c)  & = & 0,
\end{array}\right.\label{sysPR}\\
c(0,x)=c^- \in [0,1],  \quad x > 0, &\quad
&
\left\{
\begin{array}{ccl}
c(t,0) &=&c^+   \in [0,1],\\
u(t,0)&=&u^+ > 0,
\end{array}
\right.
t>0.\label{dataPR}
\end{eqnarray}
We are classically looking for a selfsimilar
solution, i.e.: $\ds c(t,x)=C(z)$,  $u(t,x)=U(z)$ with $z=\ds\frac{t}{x} >0$. The answer is given
by the three following results (\cite{BGJ07}).

\begin{proposition}
Assume for instance  that $0\leq a<c^-<c^+<b\leq 1$ and $f''>0$ in $]a,b[$. Then the only smooth
self-similar solution of (\ref{sysPR})-(\ref{dataPR}) is such that :
\begin{equation}\label{Cz}
\left \{\begin{array}{cccr}
C(z)&=&c_- ,&   0  <  z  <  z_-,\\
\ds\frac{d C}{ dz} & = &
 \ds \frac{H(C)}{z\,f''(C) },&\;   z_-  <   z  < z_+, \\
C(z) & =&   c_+, &   z_+  <  z,
\end{array}  \right.
\end{equation}
where
$ z^+=\ds\frac{H(c^+)}{u^+} $, $z^-=z^+\ds\,e^{-\Phi(c^+)}$ with $\Phi(c)=\ds\int_{c^-}^c
\frac{f''(\xi)}{H(\xi)}\,d\xi$. Moreover $u^-=\ds\frac{H(c^-)}{z^-}$  and $U$ is given by:
\begin{equation}\label{Uz}
\left \{  \begin{array}{cccr}
U(z)& = &   u_- ,&   0  <   z  <   z_- ,\\
U(z)& =& \ds\frac{H(C(z))}{z}, & \;  z_-   <   z   < z_+, \\
U(z)& =&   u_+ &   z_+  <   z.
 \end{array}  \right.
\end{equation}
\end{proposition}

\begin{proposition}
If $(c^-,c^+)$ satisfies the following admissibility condition equivalent to the Liu
entropy-condition (\cite{L76}):
$$\hbox{for all } c \hbox{ between } c^- \hbox{ and } c^+, \quad \frac{f(c^+)-f(c^-)}{c^+ - c^-}\leq
\frac{f(c)-f(c^-)}{c - c^-},$$
 then the Riemann problem (\ref{sysPR})-(\ref{dataPR}) is
solved by a shock wave defined as:
\begin{equation}\label{CUshock}
C(z)=\left\{\begin{array}{ccl}
c^- &\hbox{ if } &  0 < z < s, \\
c^+ & \hbox{ if }&  s<z
\end{array}\right.,
\qquad
U(z)=\left \{\begin{array}{ccl}
u^- &\hbox{ if }& 0< z < s,\\
u^+ &\hbox{ if } & s < z,
\end{array}  \right.
\end{equation}
 where  $u^-$ and  the speed $s$ of the shock are obtained through 
$$\frac{[f]}{u^-\,[c]} + \frac{1+h^-}{u^-}=s = \frac{[f]}{u^+\,[c]} + \frac{1+h^+}{u^+},$$
where $[c]=c^+ - c^-$, $[f]=f^+ - f^- = f(c^+) - f(c^-)$, $h^+=h(c^+)$, $h^-=h(c^-)$.
\end{proposition}

\begin{proposition}
Two states $U^-$ and $U^+$ are connected by a contact discontinuity if and only if $c^-=c^+$ (with
of course $u^-\neq u^+$), or $c^-\neq c^+$ and $f$ is affine between $c^-$ and
$c^+$.
\end{proposition}

It appears from these results that we can build a weak entropy solution of the Riemann problem
(\ref{sysPR})-(\ref{dataPR}) in a very simple way (see \cite{BGJ07}), similar to the scalar case
with flux $f$, for any data. In particular, if $f''$ has a constant sign (which is the framework in
Section \ref{sFTA}), the Riemann problem is always solved by a simple wave.


\section{ Case with  smooth  concentration}\label{skf}

System (\ref{sysad}) has the strong property that there exist  weak entropy
solutions with
{\it smooth} concentration $c$ on  $(0,T)\times(0,X)$ but not necessarily smooth
velocity $u$, for
some positive constants $T$ and $X$. Furthermore, $c$ is the solution  of a
scalar conservation law.


\subsection{Existence of weak entropy solutions with smooth
concentration}\label{sscsc}


For this section, we refer to \cite {CGM04}, \cite {CGM03}.
We have a similar result in \cite{BGJ06} but only with smooth velocity. 
 Here, we obtain by the classical  method of characteristics {\it existence and
uniqueness} of a weak
entropy solution with smooth concentration and only $L^\infty$ velocity.

\newpage
\begin{theorem}[Unique weak entropy solution 
with {\it smooth} concentration]~\\ \label{thewessc}
Let $T_0 > 0$, $ X> 0$,
$ c_0 \in W^{1,\infty}([0,X],[0,1])$,  
$ c_b \in W^{1,\infty}([0,T_0],[0,1])$, 
$ \ln u_b \in L^{\infty}([0,T_0],\R)$. 
\\
If $c_0(0)= c_b(0)$ then there exists $T\in ]0,T_0]$ such that System
(\ref{sysad})-(\ref{sysad0})  
admits a {\bf unique  } weak entropy  solution $(c,u)$ on $[0,T]\times[0,X]$ with 
\begin{eqnarray*} 
   c \in W^{1,\infty}([0,T]\times[0,X],[0,1]),
 & \qquad  
\ln u \in  L^{\infty}([0,T], W^{1,\infty}([0,X],\R)).
\end{eqnarray*}
Furthermore, for any $\psi \in C^1([0,1],\R)$, setting 
$$F'(c) =(H(c)\,G(c))^{-1} \hbox{ and } Q'=H\,\psi' + h'\,\psi,$$ 
$(c,u)$ satisfies: 
\begin{eqnarray} 
\partial_x(u\,\psi(c)) + \partial_t Q(c) = 0, 
\quad 
\partial_x(u\,G(c) ) = 0,\label{eqeisc} \\  \nonumber \\
\partial_t c + u_b(t)\,G(c_b(t))\,\partial_x F(c) = 0. \label{eqRIw}
\end{eqnarray}
\end{theorem}

\pr
we build a solution  using the  Riemann invariants and we check 
that such a solution is  an entropy solution. Next, we prove uniqueness.\\
Using the Riemann invariant $W=u\,G(c)$ ($\partial_xW= 0$) and the boundary
data we define $u$ by:
\begin{eqnarray*}
\label{numero} 
u(t,x)= \frac{u_b(t)\, G(c_b(t))}{G(c(t,x))},
\end{eqnarray*}
so $u$ is smooth with respect to  $x$.
Then, the first equation of (\ref{smooth}) can be rewritten as follows:
\begin{eqnarray} \label{eqtransportc}
\partial_t c+ \mu \,\partial_x c=0,
 &  
\quad \mbox{ with  } \quad 
&
 \ds 
   \mu = \lambda^{-1}= \frac{u}{H(c)}
         = \frac{u_b(t)\, G(c_b(t))}{H(c)\,G(c)} = \mu(t,c).
\end{eqnarray}
We solve (\ref{eqtransportc}) supplemented by initial-boundary value data
$(c_0,c_b)$ by  the standard characteristics method.  Let us define, for a
given $(\tau,x)$, $X(\cdot,\tau,x)$ as the solution of: 
\begin{eqnarray*}
\frac{dX(s,\tau,x)}{ds}=  \mu(s,c(s,X(s,\tau,x))),
 & \quad &
X(\tau,\tau,x)=x.
\end{eqnarray*} 
Since $\ds  \frac{dc}{ds}(s,X(s,\tau,x)) =0$ from (\ref{eqtransportc}), we have
\begin{eqnarray*} 
 X(s,\tau,x) =   x - b(s,\tau)\,F'(c(\tau,x))
&
\mbox{ with }  &
 b(s,\tau)=\ds \int_s^\tau u_b(\sigma) \,G(c_b(\sigma))\,d\sigma.
\end{eqnarray*}
Now, for some $T\in[0,T_0]$ defined later on, we split $\Omega= [0,T]\times[0,X]$
according to the characteristic line $\Gamma$ issuing from the corner $(0,0)$,
i.e. we define the sets  
$\Omega^\pm = \{(t,x) \in \Omega,\; \pm (x-X(t,0,0)) \geq 0 \}$. 
\\
Since $\partial_x X(t,0,x)=1- b(t,0)\,F''(c_0(x))\,\partial_x c_0(x)$, 
$b(0,0)=0$ and $ b(.,0) \in W^{1,\infty}(\Omega^+)$,
the mapping $x \mapsto X(t,0,x)$ is a Lipschitz diffeomorphism for $0\leq
t\leq T$ with $T\in]0,T_0]$ small enough. Then we define on $\Omega^+$, for each
$t\in [0,T]$, $\xi(t,x)$ such as $X(t,0,\xi(t,x))=x$.
Then we have $ c(t,x)= c_0(\xi(t,x))$ on $\Omega^+$. 
Furthermore $ \partial_t \xi = -\partial_s X/\partial_x X$ and thus $c$ is
Lipschitz continuous in time and space on $\Omega^+$.
\\
We work in a similar way on $\Omega^-$ and $c \in  W^{1,\infty}(\Omega^-)$. 
Since $c$ is continuous on $\Gamma$ from the compa\-tibility conditions $c_0(0)=c_b(0)$ 
we have $c \in  W^{1,\infty}(\Omega)$.
\\
By construction $(c,u)$ satisfies (\ref{smooth}) rewritten as follows:
\begin{eqnarray*}
\partial_x \ln u = - \partial_x g(c),
& &
 u\,\partial_x c+ H \,\partial_t c =0.
\end{eqnarray*}
These equations imply: 
$$
\partial_x u  =  - u\, \partial_x g(c) 
   = - u\, g'(c)\,\partial_x c
   =  - u \,g'(c)\, \left(-\frac{H(c)}{u}\, \partial_t c \right) 
  =  - h'(c)\, \partial_t c= - \partial_t h(c). 
$$
Now we check that $(c,u)$ satisfies  (\ref{eqeisc}).
Let $\psi$ be a $C^1$ function. Using the identity
 $ Q'=h'\psi + H \psi'$  
 and the previous equations we have:
\begin{eqnarray*} 
\partial_x(u\,\psi(c)) + \partial_t Q(c) 
 & = & 
  \psi\, \partial_x u  + u \,\psi' \,\partial_x c + Q' \,\partial_t c
= \psi \,\left( \partial_x u + h'\,\partial_t c \right) 
 +\psi'\, \left( u\partial_x c + H\,\partial_t c \right) 
\\
&= &\psi \times 0 + \psi' \times 0 = 0.
\end{eqnarray*}
By the way (\ref{eqeisc}) implies (\ref{ineqentropies}), so 
$(c,u)$ is an entropy solution of  System (\ref{sysad}).
\\
We now prove the uniqueness of such a weak entropy solution.\\
Precisely, if $ c \in W^{1,\infty}([0,T]\times[0,X],[0,1])$ and 
$\ln u \in L^\infty ((0,T),  W^{1,\infty}(0,X))$ 
satisfy  (\ref{ineqentropies}) in $\mathcal{D}'([0,T[\times[0,X[)$ with 
initial-boundary data $c_0,c_b,u_b$ then we show that $(c,u)$ is necessarily  the
previous solution built by the method of characteristics.
\\
Choosing the convex functions $\psi(c) = \pm 1$ and $\psi(c) = \pm c$ we obtain 
(\ref{sysad}).
The main ingredient to conclude the proof is the fact that $u$ admits 
 a classical partial derivative only with respect to $x$. Thus classical
computations with smooth functions  to obtain (\ref{smooth})
as in the proof of Proposition \ref{PropRi}  are still valid. 
 Now $(c,u)$ satisfies  (\ref{smooth}), which implies 
 from the beginning of the proof of Theorem \ref{thewessc}
 that $(c,u)$ is our previous solution.  
\cqfd 

\begin{remark}~
\begin{enumerate}
\item 
Notice that   $T,X$  are only  depending on $ \|\ln(u_b) \|_{L^\infty}$, $ \|
c_b \|_{W^{1,\infty}}$, $ \| c_0 \|_{W^{1,\infty}}$. Thus,
if  $\ds \left(u_b^\eps \right)_{0<\eps \leq 1}$ is a sequence of 
boundary velocity data such that $\ds \left(\ln u_b^\eps \right)$ is
uniformly bounded  in $L^\infty(0,T_0) $, and if  $(c^\eps_0), (c^\eps_b)$ are
some initial and boundary concentration data uniformly bounded in $W^{1,\infty}$
with the compatibility condition at the corner  $c_0(0)=c_b(0)$, then there
exist $T>0$ and $X>0$ and Lipschitz bounds for $c^\eps, \ln u^\eps$
on $[0,T]\times [0,X]$ independent of $\eps$.
\item 
As in \cite{BGJ06}, we have a global solution with smooth concentration 
if $\lambda$ is genuinely nonlinear (for instance  an inert case and a Langmuir
isotherm), with monotonicity assumptions on $c_0$ and $c_b$. 
 \end{enumerate}
\end{remark}


\subsection{Strong stability with respect to  velocity}\label{ssStcs}


In case of a Lipschitz continuous concentration, we now give a strong stability
result for the concentration with respect to a weak limit of the  boundary
velocity.

\begin{theorem}[Strong stability for  {\it smooth} concentration]~\\ 
\label{thess}
Let be $T_0 > 0$, $ X> 0$,
 $ c_0 \in W^{1,\infty}([0,X],[0,1])$,  
$ c_b \in W^{1,\infty}([0,T_0],[0,1])$ such that $c_0(0)= c_b(0)$, 
 and  $(\ln u_b^\eps)_{0<\eps\leq 1}$ a bounded 
 sequence in $L^\infty(0,T_0)$.
Then, there exists $T\in ]0, T_0[$
such that  System (\ref{sysad})  admits 
a  unique  weak entropy  solution 
 $(c^\eps,u^\eps)$ 
with  $ c^\eps \in W^{1,\infty}([0,T]\times[0,X],[0,1]),  
\ln u^\eps \in  L^{\infty}([0,T], W^{1,\infty}([0,X],\R))$
and 
  with  initial and boundary values:
\begin{equation} \label{sysad0eps}
\left\{   \begin{array}{ccl}
\vspace{2mm}c^\eps(0,x)&=&c_0(x) \in [0,1],  \quad x > 0,\\
\vspace{2mm}c^\eps(t,0) &=&c_b(t)   \in [0,1],\quad t>0,\\
u^\eps(t,0)&=& \ds u_b^\eps\left(t \right)>0,\quad t>0.
\end{array}\right.
\end{equation} 
If $(u_b^\eps)$ converges  towards $ \overline{u}_b$  
in $L^\infty(0,T_0)$ weak$-*$ when $\eps$ goes to $0$, 
then $(c^\eps)$ converges in \newline $L^\infty([0,T]\times[0,X])$ towards the
unique
smooth solution of 
\begin{eqnarray}  \label{eqclim}
\partial_t c + \overline{u}_b(t)\, G(c_b(t))\,  \partial_x F(c)= 0,
\qquad   c(t,0)=c_b(t), \; c(0,x)=c_0(x). 
\end{eqnarray}
Furthermore we have:
\begin{eqnarray*}
\lim_{\eps \rightarrow 0} 
\left \| u^\eps(t,x) - u_b^\eps(t)\frac{ G(c_b(t))}{G(c(t,x))}
\right\|_{L^\infty ([0,T]\times[0,X])} = 0.
 \end{eqnarray*}
\end{theorem}
\pr
thanks to Theorem \ref{thewessc}, there exists $T>0$ such that 
System (\ref{sysad}), with  initial and boundary values (\ref{sysad0eps})
admits the unique  weak entropy  solution $(c^\eps,u^\eps)$
with smooth concentration in the previous sense
on $[0,T]\times[0,X]$. \\
Since  $(c^\eps)$ is bounded in $W^{1,\infty}$, up to a subsequence,
$(c^\eps)$ converges strongly in $L^\infty $  to  $c$. 
Using  (\ref{eqRIw}) in conservative form, we can pass to the limit 
and get (\ref{eqclim}). 
Problem (\ref{eqclim}) has a unique solution by the method of characteristics.
Thus, the whole sequence $(c^\eps)$ converges. We recover the last limit 
for $u^\eps$ thanks to $\partial_x (u\,G(c)) =0$.
\cqfd

Notice that if $\overline{u}_b$ is a constant function 
for instance $u_b^\eps(t)=u_b(t/\eps)$ with $u_b$ periodic
we can compute the concentration with only using a constant velocity
(the mean velocity) as in liquid chromatography.
\\
\underline{An example from  geometric optics}:
if $u_b^\eps(t) =\ds u_b\left(t,\frac{t}{\eps} \right)$
where $u_b(t,\theta) \in L^\infty((0,T),C^0(\R/\Z))$, $\inf u_b >0$,
we have a similar result with  Equation 
 (\ref{eqclim}) for $c$
where 
$\ds 
\overline{u}_b(t)= \int_0^1 u_b(t,\theta)\, d\theta
$
and a  profile $U$:
\begin{eqnarray*}
\lim_{\eps \rightarrow 0} 
\left \| u^\eps(t,x) - U\left(t,x,\frac{t}{\eps} \right) \right\|_{L^\infty}= 0&
\mbox{ where } 
U(t,x,\theta)= \ds u_b(t,\theta) \frac{G(c_b(t))}{G(c(t,x))}.
\end{eqnarray*}


\section{Front Tracking Algorithm}\label{sFTA}


 In Section \ref{skf}, where $c$ is smooth and $ \ln u_b $ is in
 $L^\infty(0,T)$, we have seen that there exists a {\it smooth} function $v$
 such that  
  \begin{eqnarray}  \label{udecomposition}
        u(t,x) & = & u_b(t)\,v(t,x).
  \end{eqnarray}
Furthermore $c$ satisfies the {\it scalar} conservation law (\ref{eqRIw}). 
For only $BV$ data we cannot expect to obtain 
such a scalar conservation law for the concentration, except in the case of
linear isotherms. In that case, the scalar conservation law  (\ref{eqRIw}) and
System (\ref{sysad}) have the same solution for the Riemann Problem, but linear
isotherms are of a poor interest from Chemical Engineering point of view.
 The first interesting  case is the case with an inert gas and a
Langmuir isotherm, first mathematically studied in \cite{BGJ06}.

Nevertheless we guess that (\ref{udecomposition}) is still true with $v \in
BV$. From \cite{BGJ06,BGJ07} we have yet obtained $BV$ regularity with respect to
$x$ with a Godunov scheme. To get $BV$ regularity with respect to $t$ we will
use a more precise algorithm to study wave interactions, namely a Front Tracking
Algorithm (FTA).

The Front Tracking method for scalar conservation laws was introduced by 
Dafermos, \cite{D72}.  The method was extended to genuinely  nonlinear systems of 
two conservation laws by DiPerna \cite{D76}.  For our purpose, we do not use the 
generalisation to  genuinely  nonlinear systems  of any size by Bressan  \cite{Br92} or Risebro
\cite{R93}.

The  FTA is much more complicated when an eigenvalue is piecewise ge\-nu\-inely  nonlinear, see
\cite{AM07, AM04, GLF07}. Then, we restrict ourselves to the  case where $\lambda$ is
genuinely nonlinear, which allows us to treat some relevant cases from the point of view of chemical
engineering  like an inert gas with a Langmuir isotherm, two active gas with a binary Langmuir
isotherm for instance. For this purpose we work in the framework exposed in the recent and yet
classical Bressan's Book \cite{Br00}. In this framework we  assume  $f'' \geq
0$, then a Riemann
problem presents only two waves:
\begin{enumerate}
\item
a contact discontinuity with speed $0$,
\item
a rarefaction wave with speed  $\lambda > 0$
\\
or a shock wave with speed between $\lambda^-$ and $\lambda^+$, characteristic speeds associated to
the left and right states, respectively.
\end{enumerate}
Let be $\delta >0$. A $\delta$-approximate Front Tracking solution 
of  System (\ref{sysad}) is a pair  of piecewise  constant function $c^\delta(t,x),
u^\delta(t,x)$, whose jumps are located along finitely many straight lines
$t=t_\alpha(x)$ in the $t-x$ plane and approximately satisfy the entropy
conditions. For each $x >0$ and $\psi'' \geq 0$, one should thus have an
estimate of the form:
\begin{eqnarray}
\label{weakEC}
\ds \sum_\alpha 
\left( [u^\delta\,\psi(c^\delta)] - \frac{d t_\alpha}{dx}[Q(c^\delta)] \right) (t_\alpha,x)& \leq & 
{\cal O}(\delta),
\end{eqnarray}
where $[u]=u^+-u^-$ is the jump across a jump line,
and the sum is taken over all jump for $x$ fixed.
Inequality (\ref{weakEC}) implies that $(c^\delta,u^\delta)$ is ``almost an entropy 
solution'':
\begin{equation}\label{IEA}
\partial_x u^\delta \psi (c^\delta) +\partial_t \psi(c^\delta) \leq {\cal O}(\delta).
\end{equation}
That's enough to get an entropy solution ``issued from FTA'' when $\delta$ goes to zero.

Since we want to only use piecewise constant functions, it is convenient to
approximate a continuous rarefaction wave by a piecewise constant function.
For this purpose,  the rarefaction curve is dicretized with a step of order
$\delta$ and then (\ref{weakEC}) still holds.

We now briefly describe an algorithm which generates these Front Tracking approximations. The
construction starts on the initial line $x=0$ and the boundary $t=0$  by taking a piecewise constant
approximation of initial value $c_b(t), u_b(t) $ and boundary values $ c_0(x)$. 
Let $ t_1 <\cdots < t_N$, $\tilde{x}_{1} < \cdots <\tilde{x}_M$ be the points where initial-boundary
values are discontinuous. For each $\alpha = 1,\cdots,N$, the Riemann problem generated by the jump
of initial  constant values at $(t_\alpha,x=0)$  is approximately solved on a forward neighborhood
of $(t_\alpha,0)$ in the $t-x$ plane by a function  invariant on line $ t-t_\alpha= \mbox{a}\,x$,
for all positive $a$, and piecewise constant. Notice that the boundary is characteristic, then we
have only one wave associated with the speed $\lambda$ in the corner $(0,0)$.

The approximate solution $(c^\delta,u^\delta)$ can then be prolonged until $x_1 > 0$ is reached,
when the first
set of interactions between two wave-fronts takes place. If $x_1 > \tilde{x}_{1}$ we first have to
solve the characteristic boundary Riemann problem at $ (t=0,x=\tilde{x}_{1})$.
Since $(c^\delta,u^\delta)(.,x_1)$  is still a piecewise constant function, the corresponding
Riemann problems can
again be approximately solved  within the class of piecewise  constant functions. The solution  is
then continued up to a value $x_2$ where the next characteristic boundary Riemann problem  occurs or
the second set of wave interactions takes place, and so on. 

According to this algorithm, contact discontinuity fronts  travel with speed zero, shock fronts
travel exactly with Rankine-Hugoniot speed, while rarefaction fronts travel with an approximate
characteristic speed. However, one exception to this rule must be allowed if three or more fronts
meet at the same point. To avoid this situation, we must change the positive speed $\lambda$ of one
of the incoming shock fronts or rarefaction fronts. Of course this change of speed can be chosen
arbitrarily small and we have again Inequality (\ref{weakEC}).

Notice that,  for $2\times 2 $  system the number of wave-fronts  cannot approach infinity in finite
$x>0$. DiPerna   shows in \cite{D76} that the process of regenerating the solution
by solving local Riemann problems yields an approximating solution  within the class of piecewise
constant functions that is globally defined and that contains only a finite number of
discontinuities in  any compact subset of the $t-x$ quarter plane
 $ t\geq 0, x\geq 0$. We then do not consider non-physical fronts as in \cite{Br00} for general $
n\times n $ systems with $ n\geq 3$.


\section{About the shock and rarefaction curves}\label{sRP}


In this section  we state our assumptions to use the FTA
with large data. Precisely we work in classical hyperbolic case,
namely, eigenvalues are linearly degenerate or genuinely nonlinear. We  assume:
\begin{eqnarray}  \label{Hconvex}
    \lambda=\frac{H(c)}{u} 
 \mbox{ is {\bf genuinely nonlinear}}
 & \quad &
\mbox{ i.e }
  f \mbox{ is {\bf convex}  on } [0,1].  
\end{eqnarray}
Actually $\lambda$ is genuinely nonlinear for $f''\neq 0$, but since
$f=c_1q_2-c_2q_1$ (see (\ref{eqf})) we can assume that 
 $f'' > 0$ exchanging the gas labels $1$ and $2$ if necessary.
\\
Our analysis of wave interactions in Section \ref{sIE} is more precise with monotonous 
$\lambda$-wave curves, then we also assume:
 \begin{eqnarray}  \label{Hmonotonous}
 \mbox{$\lambda$-wave curves are {\bf monotonous}}.  
\end{eqnarray}
To state precisely this last assumption let us introduce some notations. Let $(c_-,L_-)$ be a left
constant state  connected to $(c_+,L_+)$ a right constant state by a $\lambda$-wave curve. In the
genuinely linear case, with Assumption (\ref{Hconvex}), $\lambda$-wave curve is  a rarefaction
curve with $ c_- < c_+$ or a shock curve with $c_-> c_+$. The sign of $[c]= c_+-c_-$ comes from the
general study of the Riemann problem in \cite{BGJ07}. From the Riemann invariant $w=\ln u + g(c) $
and  the Rankine-Hugoniot conditions a $\lambda$-wave curve can be written as follows (see
\cite{BGJ07}):
\begin{eqnarray}
    \label{eqlwave}
   [L]= L_+-L_- = \ln u_+ - \ln u_-
 & = & 
   T(c_+,c_-) = 
\left \{ \begin{array}{ll} 
       - [g]= -(g(c_+) - g(c_-)) &  \mbox{ if }  c_- < c_+ \\
          S(c_+,c_-) &  \mbox{ else } 
                \end{array} \right.
\end{eqnarray} 
We give an explicit formula for $S$ in  Lemma \ref{lRH}.
\\
 Notice that we use only one Riemann invariant, namely $c$, to write
$\lambda$-wave curves. Indeed $L=\ln u$   and $c$ have quite different behavior 
as  seen in \cite{BGJ06,BGJ07} and this paper. Furthermore we can give some simple
criterion   to  have monotonous $\lambda$-wave curves.   For instance, as
$g'=-h'/H$,  the rarefaction curve is monotonous if and only if $h$ is
monotonous. 
A chemical example, investigated in \cite{BGJ06}, is the case of an inert gas
($q_1=0$) and an active gas with a Langmuir isotherm:
   $ \ds q_2^*(c_2) = Q_2 \frac{K_2 c_2}{1+K_2 c_2}$.
For this case we have
\begin{eqnarray}
   \label{likeLangmuir}
    f'' >  0,\qquad      
    & h'<0, \qquad  &
\ds \frac{\partial S}{\partial c_-} \geq 0 
    \geq   \frac{\partial S}{\partial c_+}. 
\end{eqnarray}   

The first condition of (\ref{likeLangmuir}) gives us (\ref{Hconvex}) and the
last one gives  us (\ref{Hmonotonous}).
\\
Notice that if we exchange labels $1$ and $2$ for gas, Inequalities 
(\ref{likeLangmuir}) simply become:
$$
 f''< 0 < h', \quad
\frac{\partial S}{\partial c_-} \leq 0 
    \leq   \frac{\partial S}{\partial c_+}.
$$
Let us give some isotherm examples such that 
 (\ref{likeLangmuir}) is satisfied.
\begin{proposition} \label{PHtrue}
For the following examples, Assumptions (\ref{Hconvex}), (\ref{Hmonotonous}) are
valid: 
\begin{enumerate}
\item
one gas is inert: $q_1=0$, and the other has a concave isotherm: $q_2^{''} \leq
0$,
\item
two active gas with linear isotherms:
$\ds q_i^*(c_1,c_2)
=\ds K_i c_i,
$
$i=1,2$,
\item
 two active gas with binary Langmuir isotherms:  
$ \ds q_i^*(c_1,c_2)
=\ds \frac{Q_i K_i c_i}{1+K_1 c_1 +K_2c_2},
$
$i=1,2$,
where positive  constants
$Q_1,  Q_2 ,K_1 \geq  K_2 $ satisfy: 
$Q_1 K_1 < Q_2 K_2$.
\end{enumerate}
Furthermore, for two active gas with binary Langmuir isotherms, $\lambda$ is
genuinely nonlinear, i.e. (\ref{Hconvex}) is satisfied, if $Q_1 K_1 \neq  Q_2
K_2$.
\end{proposition}
The first case is the most classical case when only one gas is active and his
isotherm has no inflexion point, for instance the Langmuir isotherm.
\\
The second case is less interesting in chemistry and only valid when
concentrations are near  constant states. 
\\
For the third case, notice that $K_1 \geq K_2$ is not really an assumption
(exchange the labels if necessary).\\[2mm]
{\bf \textit{ Proof of Proposition} \ref{PHtrue}:} we use some technical Lemmas
postponed to Subsection \ref{teclem}. The point is to satisfy
(\ref{likeLangmuir}).

{\it 1.} \underline{Case with an inert gas:}
we have  $h=q_2, \, f(c)=-c \, h(c), \, f'=-h-ch',\, f''=-2h'-c \, h''$, which
implies $ h' =q_2' \leq 0$, $h''=q_2" \leq 0$ and then  $f''\geq 0$.
We conclude thanks to Lemmas \ref{lc+} and \ref{Lpourinertgas}.
\medskip

{\it 2.} \underline{Case with linear isotherms:}
linear isotherms are
$q_1(c)=K_1 c$, $q_2(c)=K_2(1-c)$ with $K_1 \geq 0, \, K_2 \geq 0$
then $q'_1(c)=K_1 \geq 0$, $q'_2(c)=-K_2 \leq 0$,
$h'(c)=q'_1(c)+q'_2(c)=K_1-K_2$, $f''(c)=2(K_2-K_1)$.
We assume  $K_1 \leq K_2$, then we have  $h' \leq 0 \leq f''$. Since $q_i" = 0$,
$i=1,2$, we conclude thanks to Lemmas \ref{lc+} and \ref{derivee}.
\medskip

{\it 3.} \underline{Case with a binary Langmuir isotherm:}
we have
$\ds q_1(c)=\ds \frac{Q_1 K_1 c}{D}, \, q_2(c)=\ds \frac{Q_2 K_2(1-c)}{D}$
where
$D=1+K_1 c +K_2(1-c).$
Then
$ \ds q'_1(c)=\ds \frac{Q_1 K_1(1+K_2)}{D^2} \geq 0, \, q'_2(c)=-\ds \frac{Q_2
K_2(1+K_1)}{D^2} \leq 0, $
\\ 
$h'(c)=q'_1(c)+q'_2(c) \leq 0$ if and only if $Q_1 K_1(1+K_2) \leq Q_2
K_2(1+K_1)$,
\\
$q''_1(c)=\ds \frac{2 Q_1 K_1(1+K_2)(K_2-K_1)}{D^3} \leq 0$ if and only if $K_1
\geq K_2$,
\\
$q''_2(c)=\ds \frac{2 Q_2 K_2(1+K_1)(K_1-K_2)}{D^3} \geq 0$ if and only if $K_1
\geq K_2$,
\\
$f''(c)=\ds \frac{2(Q_2 K_2-Q_1 K_1)(1+K_1)(1+K_2)}{D^3} \geq 0$ if and only if
$Q_2 K_2 \geq Q_1 K_1$.
\\
Since $Q_1K_1 \leq Q_2 K_2$, 
 we get $f"\geq 0$ and 
$\ds \frac{Q_1}{Q_2}\leq \frac{K_2}{K_1}$. 
$\ds 1 \leq \frac{1+K_1}{1+K_2}$ because $K_1\geq K_2$,
 so we have 
$\ds \frac{Q_1}{Q_2}\leq \frac{K_2}{K_1}\frac{1+K_1}{1+K_2}$, i.e. $h' \leq 0$. 
Now we conclude with Lemmas \ref{lc+} and \ref{derivee}.
\cqfd

\subsection{Technical lemmas about shock curves}\label{teclem}

We express the  shock curves as follows.
\begin{lemma}\label{lRH}
We have 
$
\exp\{S(c_+,c_-)\}  
 =\ds \frac{u_-}{u_+}
=\ds \frac{\alpha+h_-}{\alpha+h_+}, 
$ 
where
$h_\pm=h(c_\pm)$ and $\alpha=  \ds \frac{[f]}{[c]}+1$. 
\end{lemma}
\pr
first,  from the Rankine Hugoniot conditions:
$\ds \frac{[uc]}{[c+q_1(c)]}= \ds \frac{[u]}{[h]}$,
i.e. 
$[h]=\ds \frac{[u][c+q_1(c)]}{[uc]}$, we obtain 

\begin{eqnarray} \label{eqRH1}
\ds \frac{u_+}{u_-} & = &
 \ds \frac{[c+q_1(c)]-c_-[h]}{[c+q_1(c)]-c_+[h]}
\end{eqnarray}

where $[c]=c_+-c_-$ and $[h]=h(c_+)-h(c_-)=h_+-h_-$, 
and we get (\ref{eqRH1}) thanks to the following computations:

\begin{eqnarray*} 
[c+q_1(c)]-c_-[h] & = & 
  [c+q_1(c)]-c_-\ds \frac{[u][c+q_1(c)]}{[uc]}
  =\ds \frac{[c+q_1(c)]}{[uc]} \left( [uc]-c_-[u] \right)
  \\ 
& =& \ds \frac{[c] u_+}{[uc]}[c+q_1(c)], 
\end{eqnarray*}
 \begin{eqnarray*}
   [c+q_1(c)]
 - c_+ [h]  
 & = &
  [c+q_1(c)]-c_+ \ds \frac{[u][c+q_1(c)]}{[uc]}
  =\ds \frac{[c+q_1(c)]}{[uc]} \left( [uc]-c_+[u] \right)
 \\ & =& \ds \frac{[c] u_-}{[uc]}[c+q_1(c)].
\end{eqnarray*}

Rewriting (\ref{eqRH1}) we get 

\begin{eqnarray*}
 \ds \frac{u_-}{u_+} 
  & = &
\ds \frac{[c+q_1(c)]-c_+[h]}{[c+q_1(c)]-c_-[h]}
 = \ds \frac{[q_1]+[c]-c_+[h]}{[q_1]+[c]-c_-[h]}
=\ds \frac{[q_1]+[c]+c_+(h_--h_+)}{[q_1]+[c]+c_-(h_- - h_+)}
 \\
 & = & 
\ds \frac{[q_1]-c_+h_+ + [c]+h_-c_+}{[q_1]+c_-h_-+[c]-h_+c_-}
=\ds \frac{[f]+[c]+h_-[c]}{[f]+[c]+h_+ [c]}
= \ds \frac{\alpha+h_-}{\alpha+h_+},
\end{eqnarray*}
 which concludes the proof.
\cqfd

We need to know the sign of $\alpha + h_\pm$ before studying the sign of $\ds
\frac{\partial S}{\partial c_\pm}$.
\begin{lemma}\label{signe}
If $h' \leq 0$ and $ c_+ < c < c_-$ then
$\ds\alpha +h(c_+)\geq  
         \alpha +h(c)\geq \alpha +h(c_-) >0$.
\end{lemma}
\pr
since $h' \leq 0$ and $c_+ < c_-$ we have $h(c_+) \geq h(c_-)$ and it is enough
to show that 
$\ds \frac{[f]}{[c]}+1 +h(c_-) >0$.
This inequality is equivalent to $[f]+[c] +[c]h(c_-) <0$ because
$[c]=c_+-c_-<0$.
Since $f(c)=q_1(c)-ch(c)$ the inequality is equivalent to $[q_1]+[c] < c_+ [h]$.
 We know that  $q'_1 \geq 0$, 
 $ c_+< c_-$, 
$h' \leq 0$  then  $[q_1] \leq 0, \, [c]<0, \, [h] \geq 0$
 and then 
  $[q_1]+[c] < 0 <  c_+ [h]$.
\cqfd

\begin{lemma}\label{lc+}
If $h' \leq 0$, if $f$ is convex and if $c_+ < c_-$  then  we have
$\ds \frac{\partial S}{\partial c_+} (c_+,c_-)\leq 0$.
\end{lemma}
\pr
we have $S(c_+,c_-)=[L]=\ln(u_+)-\ln(u_-)=\ln(\ds \frac{u_+}{u_-})$
and 
$\ds \frac{\partial}{\partial c_+} \ds \frac{u_+}{u_-}=
\ds \frac{\partial}{\partial c_+} \ds \frac{\alpha+h_+}{\alpha+h_-}$
thanks to  Lemma \ref{lRH}.
A calculus gives
$\ds \frac{\partial}{\partial c_+}  \ds \frac{\alpha+h_+}{\alpha+h_-}
=\frac{1}{(\alpha+h_-)^2}\left( - \frac{\partial \alpha}{\partial c_+} [h]
+h'(c_+)(\alpha+h_-) \right)$. 
Now  $\ds \frac{\partial \alpha}{\partial c_+} \geq 0 $ because $f$ is convex,
next  
$[h]\geq 0$ since $ h'\leq 0$ and $ c_+ < c_-$. Lastly
 $ \alpha+h_- > 0$ from Lemma  \ref{signe} and we get $\ds \frac{\partial
S}{\partial c_+}
(c_+,c_-)\leq 0$.\cqfd
The following result concerns the case with an inert gas:
\begin{lemma}\label{Lpourinertgas}
If $q_1=0$ and $q''_2 \leq 0$ then $\ds \frac{\partial S}{\partial c_-} \geq 0$
for $c_- > c_+$.
\end{lemma}
\pr
if $q_1=0$ then $f(c)=-ch(c), \, h(c)=q_2(c)$ then $h'(c)=q'_2(c) \leq 0$.
By a direct computation and thanks to  Lemma \ref{lRH}, we have $$\ds
\frac{u_+}{u_-}=\ds \frac{[c]-c_-[h]}{[c]-c_+[h]}=\ds
\frac{[c]-c_+[h]+[c][h]}{[c]-c_+[h]}=1+ \ds \frac{1}{\ds
\frac{1}{[h]}-\frac{c_+}{[c]}}.$$
But $\ds \frac{\partial}{\partial c_-} \frac{1}{[h]} <0$, $- \ds
\frac{c_+}{[c]}$ decreases, then $\ds \frac{u_+}{u_-}$ increases with respect to
$c_-$.
\cqfd
In the  case of two active components we
have the following result:
\begin{lemma}\label{derivee}
If $ q_1^{''} \leq 0 \leq q_2^{''}$  and if $f$ is convex then 
$\ds \frac{\partial S}{\partial c_-} (c_+,c_-)\geq 0$.
\end{lemma}
\pr
let be $c$  between $c_+$ and $c_-$. From Lemma \ref{signe} we have:
\begin{eqnarray*}
  u(c)=\ds \frac{f(c_+)-f(c)}{c_+-c}+1+h(c_+)>0,
 & \quad & 
v(c)=\ds  \frac{f(c_+)-f(c)}{c_+-c} +1+h(c)>0.
\end{eqnarray*}
 We rewrite $S$ using the functions $u,v$. With Lemma \ref{lRH} we get
immediately:
 \begin{eqnarray*} 
 \ds 
 S(c_+,c_-) & =& \ln \left( \ds
\frac{[f]/[c]+1+h_+}{[f]/[c]+1+h_-}\right)\nonumber
 \\  
& =& \ln \left( \ds \frac{u(c_-)}{v(c_-)} \right).
\end{eqnarray*}
The function $f$ is convex, so $u$ is increasing.
>From equality  $f(c)=q_1(c)-ch(c)$ we have
\begin{eqnarray*} (v(c)-1)(c_+-c) & = & 
      q_1(c_+) -c_+h(c_+)-q_1(c)+ch(c) +h(c)(c_+-c) \\
  &=&q_1(c_+)-q_1(c)-c_+(h(c_+)-h(c)).
 \end{eqnarray*}
Recall that  $h(c)=q_1(c)+q_2(c)$, so we have:
\begin{eqnarray*} (v(c)-1)(c_+-c) &=& 
q_1(c_+)-q_1(c)-c_+(q_1(c_+)+q_2(c_+)-q_1(c)-q_2(c))\\
 & =&
(1-c_+)(q_1(c_+)-q_1(c))-c_+(q_2(c_+)-q_2(c)).
\end{eqnarray*}
Finally,
$v(c)-1= (1-c_+) \ds \frac{q_1(c_+)-q_1(c)}{c_+-c} -c_+ \ds
\frac{q_2(c_+)-q_2(c)}{c_+-c}$
with $0 \leq c_+ \leq 1$.
Now, $q_1$ is concave and $q_2$ is convex,
so $v$ is decreasing. 
Finally,
  $\ds \frac{u}{v}$ is increasing 
 and $\ds \frac{\partial S}{\partial c_-} \geq 0$.\cqfd

\section{Interactions estimates}\label{sIE}
%
In this section  we study  the evolution of the total variation of $L=ln(u)$,
denoted $TV L$,  through waves interactions.  It is a key point to obtain some
$BV$ bounds and  a special structure for velocity. 

Let us denote $(c_0,L_0)$, $(c_1,L_1)$, $(c_2,L_2)$, three constant states such
that:
\begin{itemize} 
\item  
 the Riemann problem with $(c_0,L_0)$ for the left state and 
  $(c_1,L_1)$ for the right state is solved by a simple wave $\mathcal{W}_1$,
\item 
 the Riemann problem  with $(c_1,L_1)$ for the left state and 
  $(c_2,L_2)$ for the right state is solved by a simple wave $\mathcal{W}_2$,
\item 
  $\mathcal{W}_1$ and  $\mathcal{W}_2$  interact.
\end{itemize}
Just after the interaction we have two outgoing waves 
  $\mathcal{W}_1^*$, $\mathcal{W}_2^*$,
 and   the intermediary constant state
  $(c_1^*,L_1^*)$. 
We denote by $TV L$ the total variation of $\ln u$ just before interaction:
$$ TV L  =  |L_0-L_1| +  |L_1-L_2|.$$
We denote by $TV L^*$ the total variation of $\ln u$ just {\it after } 
the interaction:
$$ TV L^*  =  |L_0-L_1^*| +  |L_1^*-L_2|.$$
We use similar notation for the concentration.
\\
Denote by $\alpha_-$ the negative part of $\alpha$:
  $ \alpha_-=\max(0,-\alpha)=-\min(0,\alpha) \geq 0.$
\\ 
We have the following key estimates:
\begin{theorem}[Variation on $TV \ln u$  and $TV c$
 through two waves interaction]\label{2wi}~\\
Assume (\ref{Hconvex}). Then there exists $\Gamma> 0$, a true constant such
that:
\label{thWI} 
   \begin{eqnarray}  \label{TVLquad}
 TV L^* & \leq &  TV L  + \Gamma\, |c_0 - c_1|\,|c_1 - c_2|, \\
 TV c^*  & \leq &  TV c. \label{TVcdec}
  \end{eqnarray}
 Furthermore, if (\ref{Hmonotonous}) is also satisfied then:
  \begin{eqnarray}  \label{TVLshock}
 TV L^* & \leq &  TV L  + \Gamma (c_1 - c_0)_- \times (c_2 - c_1)_-,
  \end{eqnarray}
 in addition,  if $S$, from (\ref{eqlwave}), satisfies the following triangular
inequality: 
 $$S(c_2,c_0) \leq S(c_2,c_1)+S(c_1,c_0)$$ when $ c_0> c_1> c_2$, then 
\begin{eqnarray}  \label{TVLdec}
 TV L^* & \leq &  TV L. 
\end{eqnarray}
\end{theorem}
Inequality (\ref{TVcdec}) means that the total variation of $c $ does not
increase and Inequality (\ref{TVLshock}) means that the total variation of
$\ln u $ does not increase after a wave interaction except when two shocks
interact. In this last case the increase of $TV \ln u$ is quadratic with respect
to the concentration variation. 

Such estimates are only valid when $f$ has no inflexion point. Else,
$\lambda$-wave curves are only Lipschitz and we loose the quadratic control for
the total variation of $L$.
\medskip

\underline{\bf Proof  of Inequality (\ref{TVcdec})}:
the decay of the total variation of the concentration is straightforward
since $c$ is constant  through a contact discontinuity,
 i.e. $c_1^*= c_0$ :\\ 
 $
  TV c^*=|c_2-c_1^*| + |c_1^*-c_0| 
      = |c_2-c_0 | \leq |c_0 -c_1|+|c_1-c_2| = TV c.
$ 
\cqfd 
\medskip

\underline{\bf Proof  of Inequality (\ref{TVLquad})}:
this proof  is much more complicated. We only assume (\ref{Hconvex}). The proof
is a  consequence of the following lemmas.
 
\begin{lemma} If a $\lambda$-wave 
interacts with a contact discontinuity
\label{lDDl}
then we have $TV L^* = TV L $.
\end{lemma}
\pr
it is the simplest case.  
We have  $c_1=c_2$  from the contact discontinuity, so, with $T$ defined in
(\ref{eqlwave}),
$L_1-L_0= T(c_1,c_0)=T(c_2,c_0)$ and, 
since $c_1^*= c_0$, we have
        $L_2-L_1^*=T(c_2,c_1^*)=  T(c_2,c_0)$. 
Then $$ L_2-L_1^* = L_1-L_0, $$ which implies 
$ L_2-L_1 = L_1^*-L_0$ and  
 $TV L^* = TV L .$
\cqfd

\begin{lemma} \label{TRI}
There exists a constant $\Gamma > 0$ such that, for all 
 $c_0,c_1,c_2 \in [0,1]$:
$$ |T(c_2,c_0)-T(c_2,c_1)-T(c_1,c_0)|
\mid \leq \Gamma \mid c_2-c_1 \mid \mid c_1 - c_0 |.$$
\end{lemma}
\pr
we define $R$ by 
$R(\alpha,\beta) 
  =T(c_2,c_0)-T(c_2,c_1)-T(c_1,c_0)$.
 We have to prove that 
 $R(\alpha,\beta)={\cal O}(\alpha \beta),$
where $\alpha=c_1-c_0, \, \beta =c_1-c_2$.
We denote $c=c_2, \, b=c_1, \, a=c_0$.
We have $T \in {\cal C}^3([0,1],\R)$ since $\lambda$ is genuinely nonlinear
and $T(b,b)=0$. We apply the Taylor's formula: 
\begin{eqnarray*} 
 T(c,a) & = & T(b-\beta, b + \alpha)
=T(b,b)-\beta \partial_1 T(b,b)+ \alpha \partial_2 T(b,b)
 \\  & & +
  \int_0^1(1-t)(\beta^2 \partial_1^2S + \alpha^2 \partial_2^2 T -2 \alpha \beta
\partial_{12}^2T)(b-t \beta, b+t \alpha)dt,
 \\
T(b,a) & =& T(b, b +\alpha)=T(b,b)+\alpha \partial_2T(b,b)+\int_0^1 (1-t)
\alpha^2 \partial_2^2 T(b, b+t \alpha)dt,
 \\ 
T(c,b)& =& T(b-\beta, b)=T(b,b)-\beta \partial_1 T(b,b)+\int_0^1(1-t)\beta^2
\partial_1^2 T (b-t \beta, b)dt,
 \\ 
 R(\alpha,\beta)&=&T(c,a)-T(c,b)-T(b,a)
\\ & =&
 -T(b,b)+\int_0^1(1-t) ( \beta^2 (\partial_1^2T(b-t \beta,b + t
\alpha)-\partial_1^2T(b-t \beta,b))+\\
 && \alpha^2(\partial_2^2T(b-t \beta,b+t \alpha)-\partial_2^2 T(b,b+t \alpha ))
-2 \alpha \beta \partial_1 \partial_2 T(b-t \beta,b+t \alpha) )dt.
\end{eqnarray*}
Since
\begin{eqnarray*}  \partial_1^2T(b-t \beta,b +t \alpha)-\partial_1^2T(b-t
\beta,b) & = & {\cal O}(t\alpha)={\cal O}(\alpha),\\
\partial_2^2T(b-t \beta,b+t \alpha)-\partial_2^2T(b,b+t \alpha)& = & {\cal O}(t
\beta)={\cal O}(\beta), \\
\partial_1 \partial_2 T(b-t \beta,b+t \alpha)
 & = & {\cal O}(1),
\end{eqnarray*}
we conclude that
$R(\alpha,\beta)={\cal O}(\beta^2 \alpha+\alpha^2 \beta+\alpha \beta)={\cal
O}(\alpha \beta)$.\cqfd

To conclude the proof of Inequality (\ref{TVLquad}) it suffices to use the next
lemma.

\begin{lemma}\label{llDl} If two $\lambda$-waves interact then we have 
$TV L^* \leq TVL+ \Gamma |c_{2}-c_{1}|\, |c_1-c_{0}|.$
\end{lemma}
\pr by definition of $TV L $ and $TVL^*$ it suffices  to prove that 
  $$ L_1^* = L_0 + \mathcal{O}(|c_{2}-c_{1}|\, |c_1-c_{0}|),$$
since 
$
 TV L^*= |L_2-L_1^*| + |L_1^*-L_0| 
      \leq |L_2 -L_0| +2 |L_1^*-L_0| 
       \leq TV L  + 2 |L_1^*-L_0|. 
$
\\
 Indeed, we have: 
  $L_1 -L_0 = T(c_1,c_0)$, $L_2 -L_1 = T(c_2,c_1)$, $L_2 -L_1^* = T(c_2,c_1^*)=
T(c_2,c_0) $. Next:
 $  L_2 -L_0 = T(c_2,c_1)+T(c_1,c_0)$
 and then   
 $$ L_1^* -L_0 = T(c_2,c_1)+T(c_1,c_0) - T(c_2,c_0) ,$$
 which allows us to conclude the proof of 
 Lemma \ref{llDl} with Lemma \ref{TRI}. 
\cqfd
 The proof of  Inequality (\ref{TVLquad}) is now complete. \cqfd

\underline{\bf Proof of Inequalities (\ref{TVLshock}),
                                      (\ref{TVLdec})}:
\medskip

we  assume again (\ref{Hmonotonous}) and also (\ref{likeLangmuir}) to fix the
signs. There are  more cases to study:
\begin{itemize}
\item
 first,  we have yet studied in Lemma \ref{lDDl} the  interaction 
of a shock wave or a  rarefaction wave ($\lambda$-wave) with a
  contact discontinuity  (1-wave): the contact discontinuity is ``transparent'' 
since $TV L^*= TV L$ and the concentration variation is also invariant.
 \item
 second, we  study the interaction of a shock wave with a rarefaction wave
  ($\lambda$-waves with different types): see Lemmas
 \ref{RSDR}, \ref{RSDS}, \ref{SRDR} and \ref{SRDS}. We get
  $TV L^* <  TV L$ and the concentration variation decreases. It is the only
case where $TVL$ and $TVc $ decrease.
\item 
finally, we study the interaction of two  shock waves. In this situation
  $TV L^* \geq   TV L$ and $TVc $ is invariant. 
\\
 Furthermore, if $S$ satisfies some ``triangular inequality'', we get $TV L^* =
TV L$.
\end{itemize}
In order to simplify the notations we denote by D a contact discontinuity, R a
rarefaction wave and S a shock wave. `` RD $\rightarrow$ DR "  means that a
rarefaction wave coming from the left interacts with a contact discontinuity and
produces a new left wave, namely a contact discontinuity, and a new right wave,
namely  a rarefaction.
\\
Since a contact discontinuity has a null speed and a $\lambda$-wave has a
positive speed, the only cases for $\mathcal{W}_1,\,\mathcal{W}_2$ are: RD, SD,
RS, SR and  SS.
\\
For the resulting waves   $ \mathcal{W}_1^*, \mathcal{W}_2^*$, there are 7
cases.
\\
The first two cases {RD $\rightarrow$ DR} and {SD $\rightarrow$ DS} have yet
been studied in Lemma \ref{lDDl}.

\begin{lemma}
In the case {RS $\rightarrow$ DR}\label{RSDR},
  $TVL$ decreases i.e. $TVL^* <  TVL$.
\end{lemma}
\pr
at the beginning, we have a rarefaction,
 then $c_0 < c_1$, $L_0 > L_1$, and a shock, then $c_2 < c_1$, $L_2 > L_1$.
After the interaction, we have a contact discontinuity, then $c_0=c_1^*$, and a
rarefaction, then $c_1^* < c_2$, $L_1^* > L_2$.
Finally, we have $c_0=c_1^* < c_2 < c_1$ then $g(c_0) =g(c_1^*) \leq  g(c_2)
\leq g(c_1)$. We can write  
 \begin{eqnarray*} 
 TVL &=&\mid L_0 -L_1 \mid + \mid L_1 -L_2 \mid=L_0-L_1 + L_2 - L_1,
\\
TVL^*& = & \mid L_0 -L_1^*\mid + \mid L_2 - L_1^* \mid =\mid L_0-L_1^*\mid +
L_1^*-L_2.
\end{eqnarray*}
There are two cases:
\begin{itemize}
\item
the  simplest is $L_0 > L_1^*$,
then $TVL^*=L_0-L_1^*+L_1^*-L_2=L_0-L_2 < L_0-L_1 < TVL$,
\item 
the second case is
$L_0 < L_1^*$.
Let us define $\tilde L_2$ by 
 $$ L_0-\tilde L_2=L_1^*-L_2 ,$$
then 
$ L_0-\tilde L_2=L_1^*-L_2 
=g(c_2)-g(c_1^*)=g(c_2)-g(c_0) \leq g(c_1)-g(c_0)= L_0-L_1$ because $[L]=-[g]$ 
for a rarefaction and $c_1^*=c_0$. Since shock curves are decreasing,  we know
that $\tilde L_2 > L_1$, so  
$TVL^*=L_1^*-L_0+L_1^*-L_2 =L_2-\tilde L_2 +L_0 -\tilde L_2
 <  L_2-L_1+L_0-L_1=TVL.$\cqfd
\end{itemize}

\begin{lemma}
In the case {RS $\rightarrow$ DS} \label{RSDS}  we get 
$TVL^* \leq TVL$.
\end{lemma}
\pr
this case needs the assumption $\ds \frac{\partial S}{\partial c_-} \geq 0$.
At the beginning, we have a rarefaction: $c_1 > c_0$ and $L_1  < L_0$ with a
shock: $c_2 < c_1$ and $L_2 > L_1$.
The state $(c_2,L_2)$ is connected with a shock $(c_1^*,L_1^*)$: $c_2 < c_1^*$
and $L_1^* < L_2$.
The state $(c_0,L_0)$ is connected with a contact discontinuity $(c_1^*,L_1^*)$:
$c_0 = c_1^*$.
Finally, we have  $c_2 < c_0=c_1^* < c_1.$
Then $TVL=\mid L_0 -L_1 \mid + \mid L_1 -L_2 \mid=L_0-L_1 + L_2 - L_1$
and 
$TVL^*= \mid L_0 -L_1^*\mid + \mid L_2 - L_1^* \mid =L_2-L_1^*+\mid
L_1^*-L_0\mid.$
But, with the assumption, $\ds \frac{\partial S}{\partial c_-} \geq 0$,
$S(c_2,c_0)=S(c_2,c_1^*)=L_2-L_1^* < S(c_2,c_1)=L_2-L_1$
then $L_1^* > L_1$.
\\
There are two cases:
\begin{itemize}
\item 
if $L_0 > L_1^*$ then $TV L^*=L_0-L_1^*+L_2-L_1^* < L_0-L_1+L_2-L_1=TVL$,
\item 
else
$L_0 < L_1^*$ then $TV L^*=-L_0+L_1^*+L_2-L_1^*=L_2-L_0< L_2-L_1< TVL$.\cqfd
\end{itemize}

\begin{lemma}
In the case {SR $\rightarrow$ DR}\label{SRDR} we have $TV L^* \leq TVL$.
\end{lemma}
\pr
in the beginning, we have a shock  who interacts with a rarefaction 
then $c_1 < c_0$, $L_1>L_0$ and
$c_2  > c_1$, $L_1 > L_2$.
\\
After the interaction, we have a contact discontinuity then $c_0=c_1^*$ and a
rarefaction then $c_2 > c_1^*$ and $L_1^* > L_2$.
Finally, we have $c_1 < c_0 =c_1^* < c_2.$
Since  $g' \geq 0$, we have
$g(c_1) \leq g(c_0) \leq g(c_2)$.
\\ 
For a rarefaction $[L]=-[g]$ then
$L_2-L_{1} ^{*}=g(c_1^*)-g(c_2)=g(c_0)-g(c_2)$ because $c_1^*=c_0$,
\\
$L_2-L_1=g(c_1)-g(c_2) \leq g(c_0)-g(c_2)$ because $c_1 < c_0$ and $g' \geq 0$.
\\
So we have:   $L_2-L_1 \leq  g(c_1^*)-g(c_2)=L_2-L_1^*$ and
\\
$TVL=\mid L_1-L_0 \mid + \mid L_2-L_1 \mid =L_1-L_0+L_1-L_2
\geq L_1-L_2$,
\\
$TV L^*=\mid L_1^*-L_0 \mid +\mid L_2-L_1^* \mid=\mid L_1^*-L_0 \mid
+L_1^*-L_2.$
\\
There are two cases:
\begin{itemize}
\item 
the first is $L_1^* > L_0$ then $TV L^*=L_1^*-L_0-L_2+L_1^*=2
L_1^*-L_0-L_2=-(L_2-L_1^*)+L_1^*-L_0$
$<-(L_2-L_1)+L_1^*-L_2+L_2-L_0 < -L_2+L_1-L_2+L_1+L_2-L_0 = 2L_1-L_2-L_0 =TVL$,
\item 
the second case is $L_1^* < L_0$ then $TV L^*=-L_1^*+L_0-L_2+L_1^*=L_0-L_2  \leq
L_1-L_2 \leq TVL$.
\cqfd
\end{itemize}

\begin{figure}[!ht]
\centering
\includegraphics[scale=0.5]{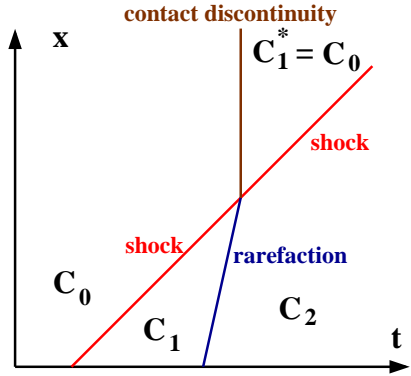}
\caption{case SR $\rightarrow$ DS. \label{SR-DS}}
\end{figure}

\begin{lemma} In the case {SR $\rightarrow$ DS}\label{SRDS},
%
 $TVL$ decreases i.e. $TV L ^* \leq TVL$.
\end{lemma}
This situation is illustrated in Fig. \ref{SR-DS}.

\begin{figure}[!ht]
\centering
\includegraphics[scale=0.4]{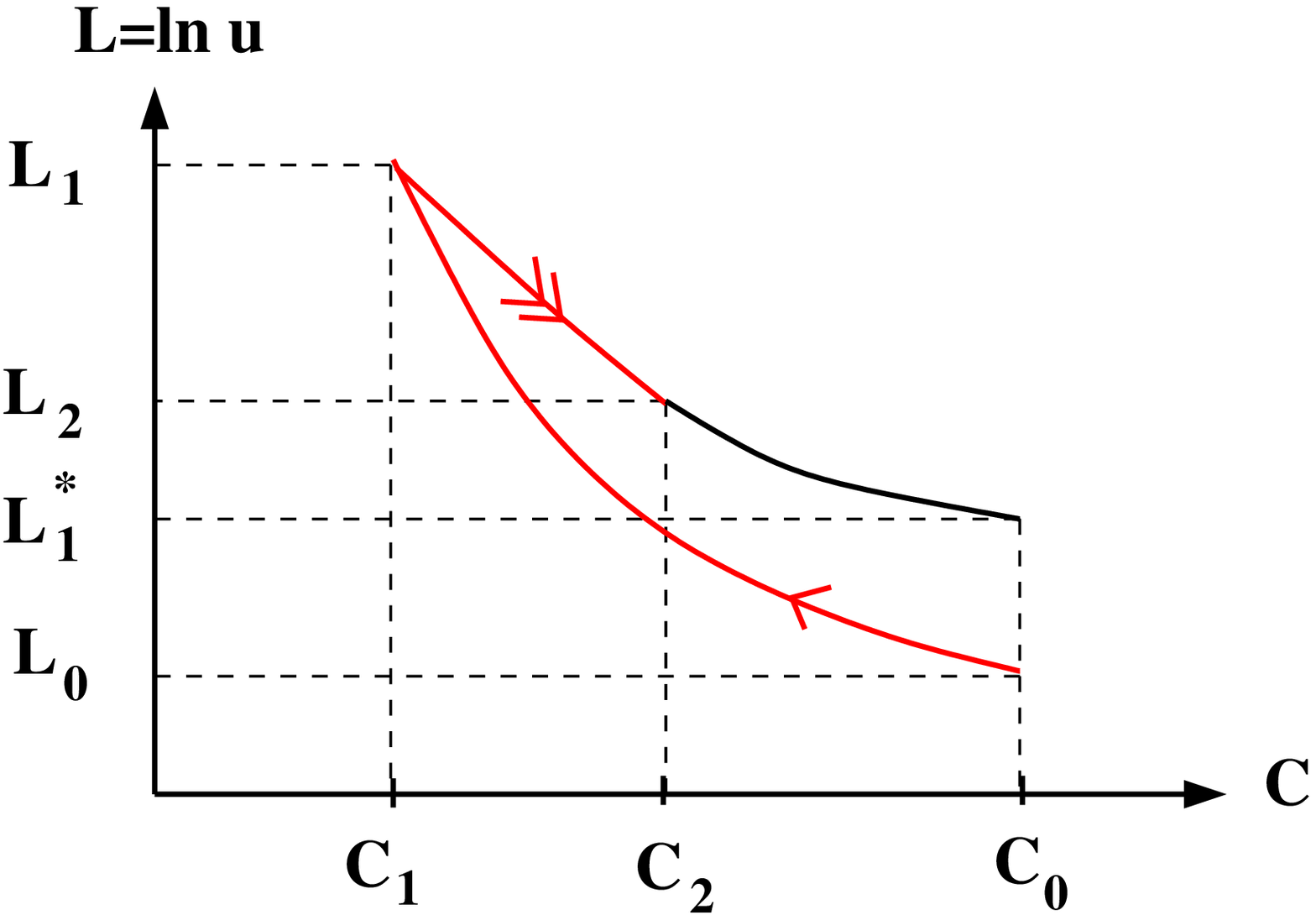}
\includegraphics[scale=0.4]{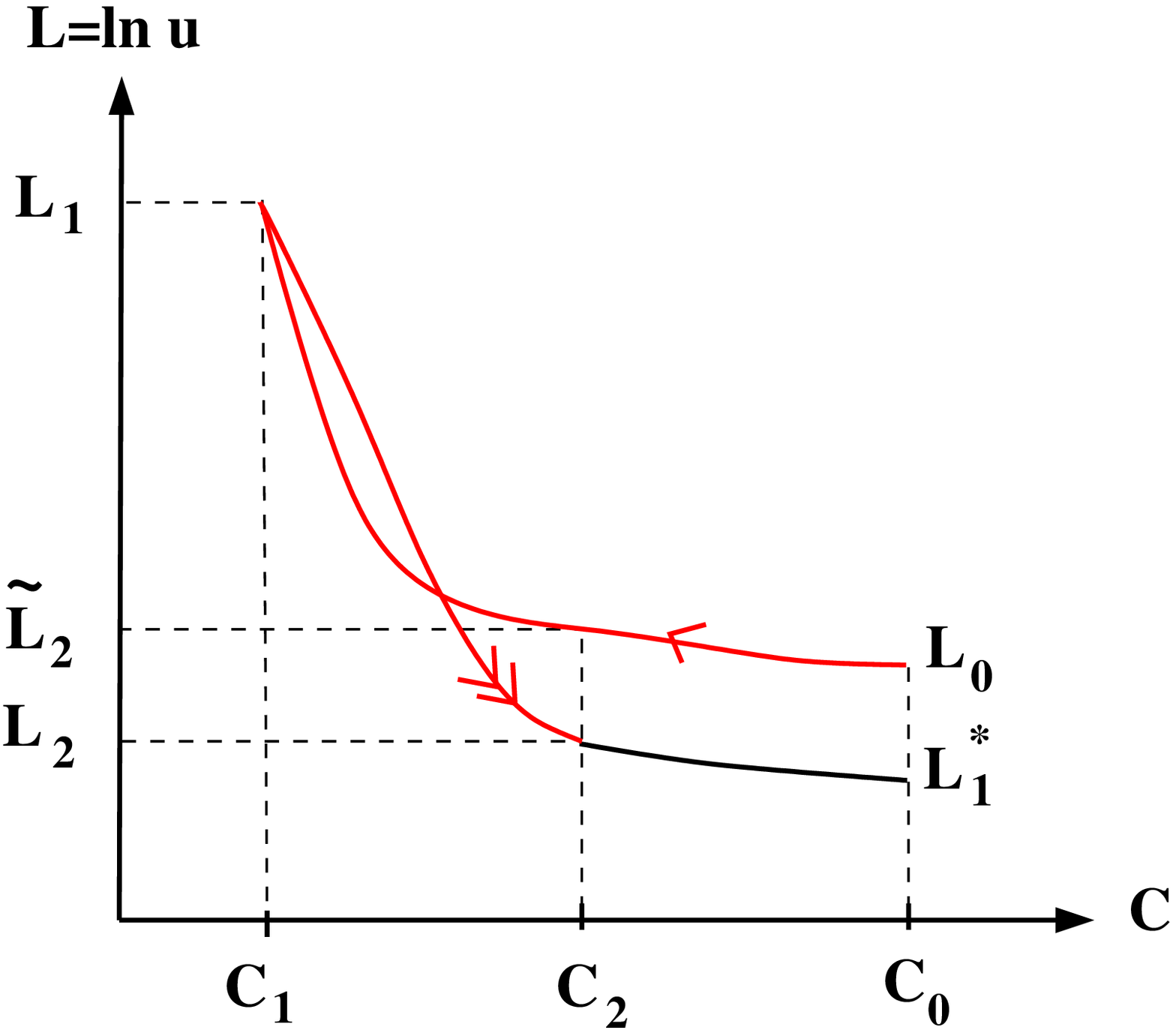}
\caption{SR $\rightarrow$ DS: first case, left, and second case,
right.\label{1and2}}
\end{figure}

\pr
it is the most difficult case.
At the beginning, we have a shock then $c_1 < c_0$ and $L_1 >L_0$.
 The shock   interacts with a rarefaction then $c_2 > c_1$ and $L_2 < L_1$.
\\
We then  have 
 $TVL =  \mid L_1 -L_0 \mid + \mid L_2-L_1 
  \mid=L_1-L_0+L_1-L_2.$
\\
The state $(c_2,L_2)$ is connected to $(c_1^*,L_1^*)$ by a shock then $c_2 <
c_1^*$ and $L_1^* < L_2$.
\\
The state $(c_0,L_0)$ is connected to $ (c_1^*,L_1^*)$ by a contact
discontinuity  then $c_0 = c_1^*$.
\\
Finally, we have $c_1 < c_2 < c_1^*=c_0$,
 $S(c_1,c_0)=S_{10}> S(c_2,c_0)=S_{20}
   =S(c_2,c_1^*)=L_2-L_1^*$,
\\
$L_1-L_0 =S_{10} > S_{20}=L_2-L_1^*$, 
because $\ds \frac{\partial S}{\partial c_+} <0$.
\\
There are two cases:
\begin{itemize}
\item 
if  $L_0 < L_1^*$ (see Fig. \ref{1and2}, left) then $ L_2 < L_1$ and
$$TV L^*= \mid L_1^* -L_0 \mid + \mid L_2-L_1^* \mid = L_1^* -L_0  +  L_2-L_1^*
=L_2-L_0< L_1-L_0 < TVL,$$
\item 
if $L_1^* < L_0$ (see Fig. \ref{1and2}, right) then we define $\tilde L_2$ by
$\tilde L_2 -L_0
=S_{20}=L_2-L_1^*< S_{10}=L_1-L_0$ and 
$TV L^* = \mid L_1^* -L_0 \mid + \mid L_2-L_1^* \mid =L_0-L_1^*+L_2-L_1^*=\tilde
L_2 -L_2+  S_{20} < L_1-L_0+L_1-L_0=TVL.$
\cqfd
\end{itemize}

\newpage
The following case is the only one where $TVL $  increases,
except if $S$ satisfies a "triangular inequality".


\begin{lemma} In the case {SS $\rightarrow$ DS}\label{SSDS}
 we have $$TV L^*= TVL+2 \max(S_{20}-S_{21}-S_{10},0) = TVL +2 \max
(L_0-L_1^*,0)\geq TVL.$$
\end{lemma}
\pr
at the beginning, we have a shock: $c_1 < c_0$ and $L_1 >L_0$. It  interacts with
an another shock:
 $c_2 < c_1$ and $L_2 > L_1$.
\\
The state $(c_2,L_2)$ is connected to $(c_1^*,L_1^*)$ by a shock then $c_2 <
c_1^*$ and $L_1^* < L_2$.
\\
The state $(c_0,L_0)$ is connected with with $(c_1^*,L_1^*)$  by a contact
discontinuity then  $c_0 = c_1^*$.
\\
Finally, we have $c_2 < c_1 < c_0 =c_1^*$
and $L_0 < L_1 < L_2.$
\\
With  $L_2-L_1=S_{21}>0$, $L_1-L_0=S_{10}>0$, 
     $L_2-L_1^*=S_{20}>0$,
we have: 
\\
$TVL= \mid L_2-L_1 \mid + \mid L_1 - L_0 \mid =L_2-L_1+L_1-L_0=S_{21}+S_{10},$
\\
$TV L^*=\mid L_2-L^*_1 \mid + \mid L_1^*-L_0 \mid
=\mid S_{20} \mid +\mid L_1^*-L_2+L_2-L_0 \mid$
$ 
=S_{20}+\mid -S_{20}+L_2-L_0\mid.$
\\
There are two cases to study:
\begin{itemize}
\item 
if $-S_{20}+L_2-L_0 \geq 0$
i.e. $S_{20}=L_2-L_1^* \leq S_{21}+S_{10}=L_2-L_0$ 
i.e. $L_0 < L_1^*$
\\
then 
$TV L^*=S_{20}-S_{20}+L_2-L_0=L_2-L_0=TVL$,
\item 
else $L_1^* < L_0$  and we have  
\begin{eqnarray*}
 TV L^* & = &S_{20}+S_{20}-L_2+L_0=2S_{20}-2L_2+2L_0+L_2-L_0
\\ & =& 2(S_{20}-(L_2-L_0))+TVL, 
  = 2(S_{20}-S_{21}-S_{10})+TVL
\\ &= &2 (L_0-L_1^*) +TVL,
\end{eqnarray*} 
\end{itemize}

which conclude the proof of Lemma \ref{SSDS}.\cqfd
The proof of Theorem \ref{thWI} is now complete.\cqfd

\section{$BV$ estimates with respect to time for the velocity} \label{suBV}
In System (\ref{un})-(\ref{deux})-(\ref{trois}), there is no partial derivative
with respect to $t$ for $u$. Nevertheless, the hyperbolicity of this
system ( with $x$ as the evolution variable) suggests that a $BV$ regularity of
the "initial" data $u_b$ for $x=0$ is propagated. Furthermore, in the case with
smooth concentration, the Riemann invariant $u\,G(c)$ suggests that when $\ln
u_b$ is only in $L^\infty(0,T)$, we can hope  $ u(t,x)/u_b(t)$ to be still $BV$
in time for almost all $x$. 
We  prove that this $BV$ structure of the velocity  is still valid with some convexity
assumptions, using a  Front Tracking Algorithm (FTA).
We  conjecture that this structure is still valid  for the general case without
convexity assumption or, better, with a piecewise genuinely nonlinear eigenvalue
$\lambda= H(c)/u$. But, in this last case, the FTA becomes very complicated 
(see Dafermos' comments in \cite{D00}).

\subsection{The case $\ln u_b \in BV(0,T)$}\label{BV}

We first precise the  notations used in the
next theorem. We define the function $c_I$ on $(0,T)$ by
\begin{eqnarray*} 
 c_I(s)= \left \{ \begin{array}{cl} 
            c_0(s)   & \mbox{  if  }0 < s < X \\
          c_b(-s)    &  \mbox{ if } 0 < -s < T
          \end{array} \right.,
\end{eqnarray*}
and we set $TV c_I = TV c_I[-T,X]$.
\\
 There exists a positive constant  $\gamma$ such that 
if $(c_-,L_-)$ is connected to $(c_+,L_+)$ by a $\lambda$-wave then
$|L_+-L_-|\leq \gamma \,|c_+-c_-|$. That is an easy consequence of 
 (\ref{eqlwave}). Indeed, it is  yet proven in \cite{BGJ06}, Lemma 3.1,  with
an inert gas, or in \cite{BGJ07}, Lemma 4.1, for two active gases.
\\
The constant $\Gamma$ comes from  Theorem  \ref{thWI}.
\begin{theorem}{\bf [Propagation of BV regularity in time for the velocity]}
           \label{thbvp}\\
Assume (\ref{Hconvex}). If $\ln u_b  \in  BV(0,T)$, if
$c_0,c_b \in BV $  and if $(u,c)$ is a weak entropy solution of 
 System (\ref{un})-(\ref{deux})-(\ref{trois}), coming from the Front Tracking
Algorithm, then $c \in BV ((0,T)\times (0,X))$ and  $ u\in
L^\infty((0,T),BV(0,X))\cap
L^\infty((0,X),BV(0,T))$. More precisely:
\begin{eqnarray*} 
     \max \left(\sup_{0 <t <T} TV_x  c(t,.)[0,X], 
         \sup_{0 <x <X} TV_t c(.,x)[0,T]  \right) 
      & \leq & TV c_I, \\
 \sup_{0 <t <T} TV_x \ln u(t,.)[0,X] & \leq & 
     TV \ln u_b +  \gamma\, TV c_I, 
                       \\
   \sup_{0 <x <X} TV_t \ln u(.,x)[0,T] & \leq & 
     TV \ln u_b 
   + 2 \gamma \,TV c_I + \frac{\Gamma}{2}\, (TV c_I)^2.
\end{eqnarray*}  
\end{theorem}
Compared to \cite{BGJ06,BGJ07}, the new result is that $u(t,x)$ is BV with respect to
 time if $u_b$ is in $BV(0,T)$ i.e. the last inequality of the Theorem \ref{thbvp}. With the 
Godunov scheme used in \cite{BGJ06,BGJ07} we do not obtain such time regularity for the velocity. It
is the reason why we use the FTA to get more precise estimates. Notice that we consider a local (in
time and space) problem for reasons of realism: we could consider a global one as well, i.e. for
$(t,x)\in(0,+\infty)^2$.
\medskip

\pr 
The easiest BV estimate on the concentration $c$  after interaction (estimate (\ref{TVcdec}) in
Theorem \ref{2wi}),  which is always valid independently of the velocity $u$, yields to a control of
$c$ in $L^\infty_tBV_x \cap L^\infty_xBV_x$  as in \cite{BGJ07}, since $\lambda$ waves  always have
a positive speed. From Lemma 4.8 of \cite{BGJ07} p.80 (ore more simply Lemma 3.1 of \cite{BGJ06} p.
557) we get $ L^\infty_{t,x} \cap L^\infty_tBV_x$ bounds for the velocity $u$. It follows, from a
natural adaptation of the estimates and compactness argument of the proof of Theorem 5.1 p 563. in
\cite{BGJ06} or Theorem 6.1 p.83 in \cite{BGJ07}, that there exists a subsequence which converges to
a solution of the initial boundary value problem with the prescribed data  $ c_0,c_b,u_b$  when
$\delta$ goes to zero, thanks to the approximate entropy inequality (\ref{IEA}). Furthermore, as in
\cite{BGJ06,BGJ07}, we recover  strong traces at $t=0$ and $x=0$.

Notice that this existence  proof  is also valid without any $BV$ assumption on the
velocity at the  boundary: we only need $\ln u_b$ in $L^\infty(0,T)$.

The $BV$  estimate with respect to time for $\ln u$, i.e. the third estimate in the theorem, is a
consequence of two following lemmas.\cqfd

Let $(u,c)$ be an entropy solution coming from FTA.
For $\delta > 0$, representing the distance from the boundary 
$x=0$ or $t=0$, let us define:
 \begin{eqnarray*} 
     L(s,\delta)   & = & \ds \left \{ \begin{array}{cc}
                          \ln u(t=|s|,x=\delta) & \mbox{if } -T < s <0 \\ 
                          \ln u(t=\delta,x=s) & \mbox{if } 0 < s <X
                     \end{array} \right., \\
  TV  L(0)  & =& \limsup_{\delta \rightarrow 0}  TV L(.,\delta)[-T,X].
 \end{eqnarray*}
 For piecewise data, $ TV L(0)$ is the total variation of $\ln u$ just before 
the first interaction. 
\begin{lemma}

 Before  wave-interactions 
 we have $TVL(0) \leq TV \ln u_b + 2 \gamma\, TV c_I$.
\end{lemma}
\pr
it suffices to prove this inequality for a piecewise constant approximate
solution issued from the FTA.
We discretize $[0,T]$ and $[0,X]$ as follows:
\\
$ T=s_1 > s_2 \cdots > s_m >s_{m+1}= 0 < s_{m+2} < \cdots  < s_N=X$.
\\
For $i=1,\cdots,m$ let us define the following piecewise approximations of $c$
and $\ln u$:
\begin{eqnarray*}
 c_{i} = \frac{1}{s_i-s_{i+1}}\int_{s_{i+1}}^{s_{i}} c_b(t)dt, 
&&
  L_{i} = \frac{1}{s_{i}-s_{i+1}}\int_{s_{i+1}}^{s_{i}} \ln (u_b(t))dt. 
\end{eqnarray*}
Since $t=0$ is a characteristic boundary we define only $c_i$
 for $i=m+1, \cdots, N-1$ by:
 \begin{eqnarray*}
 c_{i} = \frac{1}{s_{i}-s_{i+1}}\int^{s_i}_{s_{i+1}} c_0(x)dx.
\end{eqnarray*}
For $i<m$ we solve the $i^{th}$ Riemann Problem 
 with left state $(c_i,Li)$ and right state $(c_{i+1},L_{i+1})$
 and we denote by $c_i^*,L_i^*$ the intermediary state. 
Indeed $c_i^*=c_{i+1}$ since $c$ is constant through a contact discontinuity.
From Lemma 3.1 p. 557 of \cite{BGJ06} 
(or Lemma 4.1 p.78-79 of \cite{BGJ07} for two active gases) we know 
that:
\begin{eqnarray*} 
   |L_i - L_i^*| & \leq & \gamma\, |c_i -c_i^*|= \gamma \,|c_i-c_{i+1}|.
\end{eqnarray*}
We now estimate the total variation of $\ln u$ for the $i^{th}$ Riemann problem:
\begin{eqnarray*}
  |L_i - L_i^*|+ |L_i^*-L_{i+1}| & \leq &
   |L_i - L_i^*|+ \left(  | L_i^*-L_i|+ |L_i-L_{i+1}| \right)\\
 &\leq& 
 2 \gamma\,|c_i-c_{i+1}|+ |L_i-L_{i+1}|  .
\end{eqnarray*} 
Now, we look at the corner $t=0$, $x=0$ and $i=m$.
There is only a $\lambda$-wave since the boundary is characteristic.
With the left state $(c_m,L_m)$ and only $(c_{m+1})$ for the right state,
the resolution of the Riemann problem gives us a new constant value for $\ln u$,
 namely  $L_{m+1}=L^*_m$. 
 We have again the estimate 
$|L_m - L_m^*|= |L_m -L_{m+1}|  \leq   \gamma \,|c_m-c_{m+1}|$.
 So for $i=m+1, m+2,\cdots, N-1$
  we define  $L_i $  solving the characteristic Riemann problems
with the estimate:
 \begin{eqnarray*} |L_i -L_{i+1}|  \leq &  \gamma |c_i-c_{i+1}|.
\end{eqnarray*}
Summing up with respect to $i$, we obtain the total variation on $L$ just before
the first wave interaction:
\begin{eqnarray*}
 TVL & \leq & 
 \sum_{i<m} \left(  2 \gamma\,|c_i-c_{i+1}|+ |L_i-L_{i+1}|  \right) 
       +  \sum_{i\geq m}   \gamma\,|c_i-c_{i+1}|\\
 &\leq &TV \ln u_b  + 2 \gamma \,TV c_I.
\end{eqnarray*}
\cqfd
\begin{lemma}
We have the following estimate: 
$TVL \leq TV L(0)+ \ds\frac{\Gamma}{2}(TV c_I)^2.$
\end{lemma}
\pr
  we prove this estimate for any constant piecewise approximation 
 built from the FTA. The same estimate is still true passing to the limit.

First, we enumerate the absolute value of the jump  
 concentration initial-boundary value  from the left to the right:
$$\alpha_i  =  c_i - c_{i-1} \qquad i=1,\cdots,N. $$
Notice that we have $N+1$ constant states for the initial-boundary data:
  $ (c_0,L_0),\cdots, (c_N,L_N).$

From Theorem \ref{thWI},
 the increase of the total variation of $\ln u$ 
 is governed by following inequality
  $TV L^* \leq TV L + \Gamma |\alpha_{i-1}| |\alpha_i|$
 if the wave number $i-1$ interacts  with the wave number $i$.
Since $c$ is constant through a contact discontinuity
 ($c$ is a 2-Riemann invariant)
 and the jump of $c$ adds up if two $\lambda$-waves interact,
 we consider only interaction between $\lambda$-waves.
 Indeed we neglect that interaction with rarefaction  has the tendency 
to reduce $TV L$.
\\
 We measure the strength of $\lambda$-wave with the jump  of $c$ through the
wave.
We have positive or negative sign whether we have a rarefaction wave 
or a shock wave.
\\
Let $\ds \left( \alpha_i^k\right)_{1 \leq i \leq N-k}$   
be the strength of the $\lambda$-wave  number $i$ 
(labeled from the left to the right) after the interaction number $k$.
We have  $ \alpha_i^0= \alpha_i$ and denote by $j^k$ 
the index such that 
the interaction  number $k$ occurs  
 with the $\lambda$-wave number $j^k$ and $j^k+1$ where 
 $ 1 < j_k \leq N-k$.
For $1\leq i  < N-k$,
the strengths of $\lambda$-waves  after the interaction number $k > 0$
 are given by: 
\begin{eqnarray*} 
  \alpha_i^k & =& \ds 
 \left\{
       \begin{array}{cc}
          \alpha_i^{k-1} &  \mbox{ if } i < j^k  \\  
          \alpha_{i}^{k-1}+\alpha_{i+1}^{k-1} &  \mbox{ if } i = j^k  \\ 
          \alpha_{i+1}^{k-1} &  \mbox{ if } i > j^k   
       \end{array}
\right.,
\end{eqnarray*}
and the increasing of $TV L$ is less or equal than 
$\Gamma S^k$ where, from Theorem \ref{thWI}, 
$$S^0 = 0,\qquad  S^{k} = S^{k-1} + | \alpha_{i}^{k-1}||\alpha_{i+1}^{k-1}|.$$

Let us define the integers $l_i^k$ as follows:

$l^0_i = i$  and at each  interaction 
$$
l^k_i = \left\lbrace 
\begin{array}{rcl}
l^{k-1}_i&  if  & i < j^k,\\ 
l^{k-1}_{i+1}&   if  & i=j^k,..., N-k+1. 
\end{array}\right. 
$$
Notice that after each interactions with two $\lambda$-waves, there is only one  outgoing
$\lambda$-wave. Thus, the number of $\lambda$-waves decreases at each interactions, which proves
again (see
\cite{D76}) that the number of interactions is finite and the FTA is well
posed.

By induction, we see that:
$\ds \alpha_i^k = \sum_{l^k_i\leq l < l^{k}_{i+1}} \alpha_l$
where 
 $ \ds l^k_1=1 < l^k_2 < \cdots < l^k_{N-k+1}=N-k+1$,
 $l^0_i=i$ and $l^k_i$ is non decreasing with respect to $k$.
Now, from the definition of $S^k$,  we can deduce that:
\begin{eqnarray} 
  \label{eqSkind}
 \ds  S^k  & = & 
 S^{k-1} + \sum_{(i,j)\in J^k}|\alpha_i||\alpha_j|,
\end{eqnarray}
 where $ \ds J^k =\{ (i,j);\; 
     l^{k-1}_{j^k}\leq i < l^{k-1}_{j^k+1}\leq j  < l^{k-1}_{j^k+2}\}$.
\\
Let us check that: 
\begin{eqnarray}  \label{eqSk}
      S^k &=& \sum_{(i,j)\in I^k}|\alpha_i||\alpha_j|, 
\end{eqnarray} 
where 
 $ \emptyset = I^0 \subset I^1 \subset \cdots \subset
         I^{k-1} \subset I^k \subset \cdots \subset 
    I =\{(i,j);\; 1 \leq i < j \leq N\}$.

It is true for $k=0$. It is true for all $k$ if 
  $ I^{k-1} \cap  J^k = \emptyset$ and then 
 $ I^{k}= I^{k-1} \cup  J^k $.
The point is only to prove that $ I^{k-1} \cap  J^k = \emptyset$.
Terms $ |\alpha_i||\alpha_j|$ in the last sum of 
 (\ref{eqSkind}) have indexes  $i$ and $j$
 which appear in two consecutive intervals, i.e.
 $   l^{k-1}_{j^k}\leq i < l^{k-1}_{j^k+1}\leq j  < l^{k-1}_{j^k+2}$
and after, for $i=j^k$, 
  $ l^k_i = l^{k-1}_i$ and $ l^k_{i+1}=l^{k-1}_{i+2}$.
So $i$ and $j$ live in the same interval and 
 then terms $ |\alpha_i||\alpha_j|$ cannot appear again
 in $S^{k+1}$, $S^{k+2}$, \ldots,
 since  such intervals are not decreasing.
\\
The same is true for all indexes in $I^k$. 
 They can appear at most one time in $S^k$. 
We then have $ I^{k-1} \cap  J^k = \emptyset$
 and (\ref{eqSk}) is true.

 We easily estimate $S^k$, which concludes the proof:
 
 $$S^k  \leq  \sum_{(i,j)\in I}|\alpha_i||\alpha_j|
     \leq \frac{1}{2}\sum_{i=1}^N\sum_{j=1}^N |\alpha_i||\alpha_j|
       =  \frac{1}{2}\left( \sum_{i=1}^N|\alpha_i|\right)^2
      \leq  \frac{1}{2}\left( TV c_I\right)^2.$$

\cqfd

\subsection{The case $\ln u_b \in L^\infty(0,T)$}

For $\ln u_b \in L^\infty$ and $c_0,c_b \in BV$ we get a $BV$ structure
for the velocity. 

\begin{theorem}\label{thbvu}{\bf [ BV structure for the velocity]} 
We assume (\ref{Hconvex}).
\\
If $\ln u_b  \in  L^\infty(0,T)$, if $c_0, \, c_b \in BV$ and if $(c,u)$ is a weak entropy solution
issued from the FTA, then  
$$\max \left(\sup_{0 <t <T} TV_x  c(t,.)[0,X],\sup_{0 <x <X} TV_t c(.,x)[0,T] 
\right) \leq  TV c_I $$

and there exists a function $v$ and constants $\gamma,\,\Gamma>0$ such that
$u(t,x) =  u_b(t) \times v(t,x)$ with
$$ \ln v  \in \left \{ L^\infty((0,X),BV(0,T)) \cap 
L^\infty((0,T),BV(0,X))\right\}\subset BV((0,T)\times(0,X)),$$
\begin{eqnarray*}
       \sup_{0 <t <T} TV_x \ln v(t,.)[0,X] & \leq & \gamma\, TV c_I, \\ 
     \sup_{0 <x <X} TV_t \ln v(.,x)[0,T] & \leq &
                   2\,\gamma\, TV c_I + \frac{\Gamma}{2} \left(TV c_I \right)^2.
\end{eqnarray*}
\end{theorem}
The new result in this theorem is that $\ds \frac{u(t,x)}{u_b(t)} $ is BV with
respect to time, although $u_b$ is not assumed to be BV, but just in
$L^\infty$. The other regularity properties have yet been  proved in
\cite{BGJ06,BGJ07}.\\[2mm]
\pr
the first estimates  for $c$ are easily obtained as in Theorem \ref{thbvp} 
since the total variation of the concentration does not increase after an
interaction. The  existence proof of  such entropy solution follows the beginning of the
proof of Theorem \ref{thbvp} which is a natural adaptation of  the existence proof from \cite{BGJ06,
BGJ07} with only $L^\infty$ velocity.\\
We now study the new  $BV$ estimates for $v$. 
We can define $v$ by the relation $u(t,x)=u_b(t)v(t,x)$ because $u_b>0$.
Let be $ M  = \ln v$ and $M_b=\ln v(\cdot, x=0)$. The initial  total variation of $M$ on $x=0$ is
 $TV M_b =0$ since  $v(t,x=0) = 1$.
\\
We approach $u_b$ with a  piecewise constant data (thus in $BV$) 
and we show that the BV estimate for $M$ is independent of $u_b$.
Notice the fundamental relation: 
$$[L]  =  \ln u_+ -\ln u_-
       = \ln ( u_b(t)\, v_+) - \ln (u_b(t)\, v_-)
      = \ln v_+ -\ln v_-= [M].$$
The equality $[L]=[M]$ implies that the $\lambda$-waves (\ref{eqlwave}) are the
same 
 in coordinates $(c,L)$ and $(c,M)$. 
Then, Theorem \ref{thWI} is still valid replacing $L$ by $M$. 
We  then can repeat the proof of Theorem \ref{thbvp} 
to get $BV$ estimates for $v$.
\cqfd


\section{Weak limit for velocity with  $BV$ concentration}\label{sks}


When $c$ is only in $BV$, we cannot reduce System (\ref{sysad}) 
to a scalar conservation law for $c$ as in section \ref{skf}. 
Indeed, since the shock speeds depend on the velocity, we have a true $
2\times 2$ hyperbolic system. 
 Nevertheless we can state following stability result.

\begin{theorem}[Stability with respect to weak limit for the velocity 
 \label{thsswlv}]~\\
Let  $ \left( \ln (u_b^\eps) \right)_{0<\eps<1}$ be a bounded sequence 
 in $L^\infty(0,T)$, such that 
$$u_b^\eps     \rightharpoonup  \overline{u}_b 
  \mbox{ in } L^\infty(0,T) \mbox{ weak *}.$$
Let be $c_0\in BV((0,X),[0,1])$ and  $c_b \in BV((0,T),[0,1])$.
Let $(c^\eps,u^\eps)$ be a weak entropy solution of  System (\ref{sysad})
on $(0,T)\times (0,X)$ issuing from  the FTA
with  initial and boundary values:
\begin{equation*} \label{sysad0epsBV}
\left\{   \begin{array}{ccl}
\vspace{2mm}c^\eps(0,x)&=&c_0(x),  \quad X> x > 0,\\
\vspace{2mm}c^\eps(t,0) &=&c_b(t),\quad T> t>0,\\
u^\eps(t,0)&=&u_b^\eps\left(t \right),\quad T>t>0.
\end{array}\right.
\end{equation*}
Then, there exists $ (u(t,x), c(t,x))$, weak entropy solution of  System
(\ref{sysad})
supplemented by initial and boundary values:
\begin{equation*} \label{eqcu0}
\left\{   \begin{array}{ccl}
\vspace{2mm}c(0,x)&=&c_0(x),  \quad x > 0,\\
\vspace{2mm}c(t,0) &=&c_b(t),\quad t>0,\\
u(t,0)&=&\overline{u}_b(t),\quad t>0,
\end{array}\right.
\end{equation*}
such that, when $\eps$ goes to $0$ and  up to a subsequence:
\begin{eqnarray*} 
  c^\eps(t,x)  & \rightarrow & c(t,x) 
      \mbox{ strongly in } L^1(  [0,T]\times[0,X]), \\ 
 u^\eps(t,x)  & \rightharpoonup & u(t,x) 
      \mbox{ weakly in } L^\infty(  [0,T]\times[0,X]) \mbox{ weak *}, \\
 u^\eps(t,x)  &  =  & u_b^\eps(t)\times  v(t,x) + o(1) 
      \mbox{ strongly in } L^1(  [0,T]\times[0,X]), 
      \mbox{ where } \ds v(t,x)= \frac{u(t,x)}{\overline{u}_b(t)}.
\end{eqnarray*} 
\end{theorem}
For the convergence of the whole sequence we need  the  uniqueness 
 of the entropy solution for initial-boundary value problem:
(\ref{sysad}), (\ref{sysad0}).
\\
\pr
from Theorem \ref{thbvu} we know that $u^\eps (t,x) = u_b^\eps(t)
v^\eps(t,x)$ where the sequences $(\ln v^\eps )_{0<\eps}$ 
 and  $( c^\eps )_{0<\eps}$ are uniformly bounded  in $BV((0,T)\times(0,X))$.
 Then, up to a subsequence, we have the following strong convergence 
 in $ L^1((0,T)\times(0,X))$: 
  $ v^\eps  \rightarrow  v $,
 $  c^\eps  \rightarrow  c $. 

$(c^\eps,u^\eps)$ is a weak entropy solution for (\ref{sysad}) means
for all $\psi$ such that $\psi"\geq 0$ and $Q$ such that 
 $ Q'=h'\psi+ H\psi'$ we have in distribution sense:
$ 
\partial_x\left(u^\eps(t,x) \psi(c^\eps) \right) 
     +\partial_t Q(c^\eps)   \leq   0,
$
which is rewritten as follows:
$
\partial_x\left(u_b^\eps(t)v^\eps(t,x) \psi(c^\eps) \right) 
     +\partial_t Q(c^\eps)   \leq   0.
$ Passing again to the  weak-limit  against a strong limit we get:
$
\partial_x\left(\overline{u}_b(t)v(t,x) \psi(c) \right) 
     +\partial_t Q(c)   \leq   0.$ 
i.e.  $(c,u=\overline{u}_b\times v)$ is a weak entropy solution for System 
 (\ref{sysad}). 
We also can pass to the limit on initial-boundary data.

Since there exists $\delta $ such that $ 0 < \delta  < u_b^\eps < \delta^{-1}$, 
$ v^\eps(t,x) \rightarrow v(t,x)$ means 
 $ u^\eps(t,x)/u_b^\eps(t) -  v(t,x) \rightarrow 0$ 
and also means 
 $ u^\eps(t,x)- u_b^\eps(t) \times  v(t,x) \rightarrow 0$,  
which concludes the proof.
\cqfd

\underline{ An example of high oscillations 
 for velocity}:
as an example of weak limit we  consider the case of high oscillations 
 for velocity on the boundary.

Let be $u_b(t,\theta) \in L^\infty((0,T),C^0(\R/\Z,\R))$,
$\overline{u}_b(t)= \ds \int_0^1 u_b(t,\theta)d\theta$ and assume 
$\inf u_b > 0$. 
 With  $u_b^\eps(t) = \ds u_b\left(t,\ds \frac{t}{\eps} \right)$
and the same notations as in Theorem \ref{thsswlv} we have: 
\begin{itemize}
\item 
 first, oscillations do not affect the behavior  of the concentration 
 since $(c^\eps) $ converges strongly in $L^1 $ towards $c$ and 
 the limiting system depends only on the average $\overline{u}_b$ and not on
oscillations;
\item 
second, $(u^\eps)$ converges weakly towards 
$\overline{u}_b(t) \times v(t,x)$
 and we have a strong profile for $u^\eps$:
$$\lim_{\eps \rightarrow 0} 
\left \| u^\eps(t,x) - U\left(t,x,\frac{t}{\eps} \right)
\right\|_{L^1((0,T)\times(0,X))} = 0,
\mbox{ where } 
U(t,x,\theta)= \ds u_b(t,\theta) \times v(t,x).$$
\end{itemize}



\end{document}